\documentclass[11pt, letterpaper, reqno]{amsart}
\usepackage{graphicx, amssymb, color}
\usepackage{amsmath}
\usepackage{latexsym}
\usepackage{hyperref}
\usepackage{dsfont}
\usepackage{amsthm}

\numberwithin{equation}{section}

\def\be{\begin{eqnarray}}
\def\ee{\end{eqnarray}}

\def\b*{\begin{eqnarray*}}
\def\e*{\end{eqnarray*}}

\newtheorem{Theorem}{Theorem}[part]
\newtheorem{Definition}{Definition}[part]
\newtheorem{Proposition}{Proposition}[part]

\newtheorem{Lemma}{Lemma}[part]
\newtheorem{Corollary}{Corollary}[part]
\newtheorem{Remark}{Remark}[part]
\newtheorem{Example}{Example}[part]

\makeatletter \@addtoreset{equation}{section}

\@addtoreset{Definition}{section}

\@addtoreset{Theorem}{section}

\@addtoreset{Proposition}{section}

\@addtoreset{Property}{section}

\@addtoreset{Assumption}{section}

\@addtoreset{Corollary}{section}

\@addtoreset{Lemma}{section}

\@addtoreset{Remark}{section}

\@addtoreset{Example}{section}

\def \C{\mathbb{C}}

\def \E{\mathbb{E}}

\def \M{\mathbb{M}}
\def \N{\mathbb{N}}
\def \P{\mathbb{P}}
\def \Q{\mathbb{Q}}
\def \R{\mathbb{R}}
\def \S{\mathbb{S}}

\def\Ac{\mathcal{A}}
\def\Bc{\mathcal{B}}
\def\Cc{\mathcal{C}}
\def\Dc{\mathcal{D}}

\def\Ic{\mathcal{I}}
\def\Jc{\mathcal{J}}

\def\Lc{\mathcal{L}}
\def\Mc{\mathcal{M}}
\def\Nc{\mathcal{N}}

\def\Rc{\mathcal{R}}
\def\Sc{\mathcal{S}}

\def\vo{\overline{v}}

\def\Xh{{\hat X}}

\def\Jt{{\tilde J}}

\def\Nt{{\tilde N}}

\def\FF{\mathbf{F}}

\def\TT{\mathbf{T}}

\def\Tr#1{{\rm Tr}\left[#1\right]}

\def \Frac{\displaystyle\frac}

\def\no{\noindent}

\def\={\;=\;}
\def\.{\;.}

\def\eps{\varepsilon}

\def\reff#1{{\rmfamily(\ref{#1})}}

\def \ep{\hbox{ }\hfill{ ${\cal t}$~\hspace{-5.1mm}~${\cal u}$   } }

\def \proof{{\noindent \bf Proof. }}
\def \ep{\hbox{ }\hfill$\Box$}

 \def\normeL2#1{\left\|{#1}\right\|_{L^2}}
\newcommand {\lb} {\lambda}
\newcommand {\Chi} {{\bf \raise 1.5pt \hbox{$\delta$}}}

\newcommand{\nb}{\nabla}
\def\1{\mathds{1}}

\def\what{\widehat}

 \date{February 11, 2011}

\begin{document}
\title[A Probabilistic Scheme for Fully Nonlinear Nonlocal Parabolic PDEs]{A Probabilistic Scheme for Fully Non-linear Non-local Parabolic PDEs with singular L\'evy measures}
\author[]{Arash Fahim}
\address[Arash Fahim]{CMAP,, Ecole Polytechnique, Paris}
\email{fahimara@umich.edu}
\thanks{I am grateful to thank Nizar Touzi and Xavier Warin for their fruitful comments and suggestions.
}


 \maketitle

\begin{abstract}
We introduce a Monte Carlo scheme for fully nonlinear parabolic nonlocal PDE's whose nonlinearity in of Hamilton-Jacobi-Bellman-Isaacs (HJBI for short). We avoid the difficulties of infinite L\'evy measure by truncation of the L\'evy integral. The first result provides the convergence of the scheme for general parabolic nonlinearities. The second result provides bounds on the rate of convergence for  concave (or equivalently convex) nonlinearities. For both results, it is crucial to choose truncation of the infinite L\'evy measure  appropriately dependent on the time discretization. We also introduce a Monte Carlo Quadrature method to approximate the nonlocal term in the HJBI nonlinearity.
\end{abstract}

\section{Introduction} 
\label{Sec:intro}
Nonlocal fully nonlinear PDE's arise from stochastic optimization problems for controlled jump-diffusion processes e.g.  problem of portfolio optimization in L\'evy markets. There are only few examples with explicit and quasi-explicit solution, e.g. \cite{benthkarlsenreikvam1} and \cite{benthkarlsenreikvam}. In some applications the dimension of the PDE is so large that classical algorithms such as finite difference fail to approximate the solution in a reasonable time. The advantage of Monte Carlo methods is that they are less sensitive with respect to dimension on comparison with other  methods. 

The present paper generalizes the probabilistic numerical approximation scheme in \cite{ftw} for fully nonlinear parabolic nonlocal PDE's. The nonlocal PDE's (sometimes refer to as integro-partial differential equations) we study are of Hamilton-Jacobi-Bellman-Isaacs  type (HJBI for short). In  \cite{ftw} a Monte Carlo scheme is presented for parabolic fully nonlinear local PDE's and for a class of nonlinearities  convergence and a rate of convergence are established. Moreover, it is shown that the error due to estimation of the conditional expectations does not change the results if enough number of sample paths are used. 

The assumption we make is analogous to \cite{ftw}, i.e. one can separate the equation into a non-degenerate elliptic linear operator and a degenerate elliptic fully nonlinear operator and the diffusion coefficient in the linear operator dominates the partial gradient of the nonlinearity with respect to its Hessian component (See assumption {\bf F (ii)}). Then, we use the time discretization of a the jump-diffusion process corresponding to the linear part to approximate the derivatives and integral term in the HJBI nonlinearity. The separation into linear and nonlinear part is arbitrary up to the satisfaction of assumptions {\bf F}. Moreover, we assume that the IPDE satisfies a comparison principle in a suitable class of functions, which assures the uniqueness of the solution in the same class.

In this paper, we use a Monte Carlo approximation of the integral with respect to L\'evy measure which appears in the HJBI nonlinearity, which in this paper is referred to as {\it Monte Carlo Quadrature (for simplicity MCQ)}. After truncation of infinite L\'evy measure near zero, we treat the jumps using the idea in \cite{be} for finite activity jump-diffusion processes. We introduce bounds for the truncation error with respect to the derivatives of integrand and truncation level.

Although MCQ is independent of the numerical scheme, we choose to approximate the L\'evy integral inside the nonlinearity by MCQ. In this case, we also need to choose appropriate truncation bound with respect to time step which retains the convergence and rate of convergence as in the local case in \cite{ftw}.

The idea of the proof is captured from \cite{barlessouganidis} and \cite{ftw} for the convergence result and from \cite{biswasjakobsenkarlsen1} and \cite{ftw} for the rate of convergence. If the L\'evy measure is finite, the convergence result is obtained in \cite{fahim0}. However, the same result can not be easily generalized to the infinite L\'evy measure, where we need to conquer the new difficulties due to lack of Lipschitz continuity of nonlinearity.

Our first result concerns the convergence of the approximate solution obtained from the Monte Carlo scheme to the viscosity solution of the final value problem. We show that if we truncate the infinite L\'evy measure near zero such that the truncation level approaches fast enough to zero as the time discretization step goes to zero, then the approximate solution converges to a function which satisfies the equation in viscosity sense, which by comparison principle is unique.

The second result provides  the rate of convergence in the case of concave nonlinearity. The proof of the rate of convergence uses the results in \cite{biswasjakobsenkarlsen1} and \cite{biswasjakobsenkarlsen2} which generalizes the result of \cite{barlessouganidis} to nonlocal HJBI PDE's. While the assumptions we made verifies the existence and uniqueness of the solution of HJBI equation,  the solution of the nonlocal PDE  is approximated by a regular sub solution and  a super solutions. Plugging the regular sub solution and super solution into scheme and then by the usage of consistency estimate, we provide an upper and a lower bound for the rate of convergence. Here, we also need to impose a stronger condition on the truncation of the l\'evy measure than in the convergence result.

Finally, as mentioned in \cite{ftw} for nonlocal case, it is worth noticing the relation with the generalization of \cite{ks} to nonlocal case introduced in  \cite{imbertserfaty} which provides a deterministic game theoretic interpretation for fully nonlinear parabolic problems. The game consists of two players. At each time step in a predetermined  time horizon, one tries to maximize her gain and the other to minimize it by imposing a penalty term to her gain. More precisely, she starts in an initial position and chooses a vector $p$, a matrix $\Gamma$,   and  a function $\varphi$. Then, he will plug an arbitrary vector $w$  together with $p$, $\Gamma$ and $\varphi$ in a nonlinear penalty term which she should pay to be allowed to change her position by taking one step with appropriate length in the direction of vector $w$. At the final stage, she will earn as much as a function of her final position. As time step goes to zero, her value function at any time  and any position will converge to the solution of a fully nonlinear parabolic IPDE whose nonlinearity relates to the penalty term. Vector $p$, a matrix $\Gamma$ and a function $\varphi$  represent the first and second derivatives and the solution function, respectively.

The paper is organized as follows: In Section \ref{Sec:features}, the problematic features of nonlocal fully nonlinear PDE is discussed on a na\"ive generalization of the Monte Carlo method from local case in \cite{ftw} to nonlocal case. In Section \ref{Sec:mcq} the Monte Carlo quadrature (MCQ) is presented as a purely Monte Carlo approximation of  with the error analysis. Section \ref{Sec:asymp} contains the results of convergence and asymptotic properties of the scheme.

\no {\bf Notations}\quad 
For scalars $a,b\in\R$, we write $a\wedge b:=\min\{a,b\}$, $a\vee b:=\max\{a,b\}$, $a^-:=\max\{-a,0\}$, and $a^+:=\max\{a,0\}$.\\
By $\M(n,d)$, we denote the collection of all $n\times d$ matrices with real entries. The collection of all symmetric matrices of size $d$ is denoted $\S_d$, and its subset of nonnegative symmetric matrices is denoted by $\S_d^+$.
For a matrix $A\in\M(n,d)$, we denote by $A^{\rm T}$ its transpose. For $A,B\in\M(n,d)$, we denote $A\cdot B:={\rm Tr}[A^{\rm T}B]$. In particular, for $d=1$, $A$ and $B$ are vectors of $\R^n$ and $A\cdot B$ reduces to the Euclidean scalar product.\\
 We denote by $\Cc_d$, the space of bounded continuous functions from $[0,T]$ to $\R^d$.
For a suitably smooth function $\varphi$ on $Q_T:=(0,T]\times\R^d$, we define 
$|\varphi|_\infty:=\sup_{(t,x)\in Q_T}|\varphi(t,x)|$ and $|\varphi|_1:=|\varphi|_\infty+\mathop{\sup}\limits_{Q_T\times Q_T}\frac{|\varphi(t,x)-\varphi(t',x')|}{(x-x')+|t-t'|^\frac{1}{2}}$.

\section{Preliminaries and features for nonlocal PDE's}
\label{Sec:features}
Let $\mu$, $\sigma$ be functions from $[0,T]\times\R^d$ to $\R^d$ and $\M(d,d)$ and $\eta$ be a function from $[0,T]\times\R^d\times\R^d$ to $\R^d$ and $a=\sigma^\text{T}\sigma$.
Suppose the following nonlocal Cauchy problem:
\be
 \label{equationnl}
 &&-\Lc^X v(t,x)- F\left(t,x,v(t,x),D v(t,x),D^2v(t,x),v(t,\cdot)\right)= 0,~~\mbox{on}~[0,T)\times\R^d,~~~~~~~~~~~~~\\
&&v(T,\cdot)=g, ~~~~~~~~~~~~~~~~~~~~~~~~~~~~~~~~~~~~~~~~~~~~~~~~~~~~~~~~~~~~~~~~~\mbox{on}~\in\R^d. \label{terminalnl}
\ee
where 
$F:\R_+\times\R^d\times\R\times\R^d\times\S_d\times\Cc_d\rightarrow\R$ and $\Lc^X$ given by:
\b*
 \Lc^X\varphi(t,x)
 &:=&
 \biggl(\frac{\partial\varphi}{\partial t}+\mu\cdot D\varphi+\frac{1}{2}a\cdot D^2\varphi\biggr)(t,x)\\
 &&+\int_{\R^d_*}\biggl(\varphi(t,x+\eta(t,x,z))-\varphi(t,x)-\mathds{1}_{\{|z|\le1\}}D\varphi(t,x)\eta(t,x,z)\biggr)d\nu(z).
\e*
$\Lc^X$ is the infinitesimal generator of a jump-diffusion, $X_t$, satisfying SDE:
\b*
dX_t=\mu(t,X_t)dt+\sigma(t,X_t)dW_t+\int_{\{|z|>1\}}\!\!\!\!\!\!\!\!\eta(t,X_{t-},z)J(dt,dz)+\int_{\{|z|\le 1\}}\!\!\!\!\!\!\!\!\eta(t,X_{t-},z)\Jt(dt,dz),
\e*
where $J$ and $\tilde J$ are respectively a Poisson jump measure and compensated Poisson jump measure describing the jumps of process $X_t$ and are related to $\nu$ through
\b*
\nu(A)
&=&
\E\left [\int_AJ([0,1],dz)\right ]\\
\Jt(dt,dz)
&=&
J(dt,dz)-dt\times \nu(dz).
\e*
For more detailes on jump-diffusion processes, see \cite{bass} and the references therein or the classic work of \cite{strook}.

The classical solution for the problem \reff{equationnl}-\reff{terminalnl} does not exist in general and therefore we appeal to the notion of viscosity solutions for nonlocal parabolic PDE's. 
\begin{Definition}
$\bullet$ The viscosity sub(super)-solution of \reff{equationnl}-\reff{terminalnl} is a upper semi-continuous (lower semi-continuous) function $\underline v$($\overline v$)$:[0,T]\times\R^d\to\R$ such that:
\begin{enumerate}
 \item for any $(t_0,x_0)\in[0,T)\times\R^d$ and any smooth function $\varphi$ with:
\b*
0=\mathop{\max}(\mathop{\min})\{\underline v-\varphi\}=(\overline v-\varphi)(t_0,x_0)
\e*
We  have:
\b*
0&\ge(\le)&-\Lc^X\varphi(t_0,x_0)-F\left(\cdot,\varphi,D \varphi,D^2\varphi ,\varphi(\cdot)\right)(t_0,x_0).
\e*
\item $g(\cdot)\ge\underline v(T,\cdot)(\le\overline v(T,\cdot))$.
\end{enumerate}
 The  function $v$ which is both viscosity sub and super solution, is called viscosity solution of \reff{equationnl}-\reff{terminalnl}.\\
$\bullet$ We say that \reff{equationnl} has comparison for bounded functions if for any bounded upper semi-continuous viscosity super-solution $\vo$ and any bounded lower semi-continuous sub-solution $\underline v$, satisfying 
\b*
\vo(T,\cdot)\ge \underline v(T,\cdot),
\e* 
we  have $\vo\ge \underline v$ on $[0,T]\times\R^d$.
\end{Definition}

\subsection{Discretization of the jump-diffusion process}
\label{Sec:discrete}
Our purpose is to introduce a Monte Carlo method which approximates the solution of problem \reff{equationnl}-\reff{terminalnl}. For this purpose, we first need to provide a discretization for the process $X$. 

Suppose that $h=\frac{T}{n}$, $t_i=ih$, and $\kappa\ge0$. We define the Euler discretization of jump-diffusion process $X_t$  with truncated L\'evy measure by:
\be\label{euler}
&&\Xh_h^{t,x,\kappa}=x+\tilde\mu(t,x)h+\sigma(t,x)W_h+\int_{\{|z|>\kappa\}}\eta(t,x,z)\Jt([0,h],dz),\\
&&\Xh_{t_{i+1}}^{x,\kappa}=\Xh_h^{t_i,\Xh_{t_i}^{x,\kappa},\kappa}~~\mbox{and}~~\Xh_0^{x,\kappa}=x.
\ee
where $\tilde\mu(t,x)=\mu(t,x)+\int_{\{|z|>1\}}\eta(t,x,z)\nu(dz)$ and we make the choice of  $\kappa=0$ when $\nu$ is a finite measure.
Let $\Nt_t^\kappa$ and $N_t^\kappa$ be respectively the Poisson process derived from jump measure $J$ by counting all jumps of size greater than $\kappa$ which happen in time interval $[0,t]$ and its compensation, i.e.
\be\label{poisson}
N_t^\kappa=\int_{\{|z|>\kappa\}}J([0,t],dz)&\text{and}&\Nt_t^\kappa=\int_{\{|z|>\kappa\}}\Jt([0,T],dz).
\ee
One can write the jump part of $\Xh_h^{t,x,\kappa}$ as a compound Poisson process (see for example \cite{conttankov})
\be\label{compoundpoisson}
\Xh_h^{t,x,\kappa}&=&x+\mu_\kappa(t,x)h+\sigma(t,x) W_h+\sum_{i=1}^{N_h^\kappa}\eta(t,x,Z_i),
\ee
where $\mu_\kappa(t,x)=\mu(t,x)-\int_{\{\kappa<|z|\le1\}}\eta(t,x,z)\nu(dz)$, $Z_i$s are i.i.d. $\R^d_*-$valued random variables, independent of $W$ and $N^\kappa$, and distributed as $\1_{\{|z|>\kappa\}}\frac{1}{\lb_\kappa}\nu(dz)$.

\subsection{The scheme for nonlocal fully nonlinear parabolic PDE's}
\label{Subsec:problems}
In this section, we introduce a probabilistic scheme by following directly the same idea as the scheme for the local PDE's. Then, we consider some problems which prevents us to utilize the scheme in many interesting applications. Therefore, we introduce a modified version of the scheme which works for the class of nonlinearities of HJB type (Hamilton-Jacobi-Bellman).

Following the same idea as in \cite{ftw}, one can obtain the following immature scheme.
 \be\label{immaturescheme}
 v^{h}(T,.)=g
 &\mbox{and}&
 v^{h}(t_i,x)=\TT_{h}[v^{h}](t_i,x),
 \ee
where for every function $\psi:\R_+\times\R^d\longrightarrow\R$ with exponential growth:
 \be\label{TT}
 \TT_{h}[\psi](t,x)\!\!\!\!\!&:=&\!\!\!\!\!\E\left[\psi\left(t+h,\hat X_h^{t,x}\right)\right]
                   +hF_{h}\left(t,x,\Dc_h\psi,\psi(t+h,\cdot)\right),\\
\nonumber&&~\Dc_h\psi:=\left(\Dc_h^0\psi,\Dc_h^1\psi,\Dc_h^2\psi\right),
 \ee
where \be\label{hermit}
 \Dc_h^k\psi(t,x)
 &:=&
 \E\left[\psi(t+h,\hat X^{t,x,\kappa}_h)H^h_k(t,x)\right], k=0,1,2,
 \ee
 where
  \b*
 H^h_0=1,
 &H^h_1=\left({\sigma^{\rm T}}\right)^{-1}\;\Frac{W_h}{h},&
 H^h_2=\left({\sigma^{\rm T}}\right)^{-1}\;\frac{W_hW^{\rm T}_h-h\mathbf{I}_d}{h^2}\;\sigma^{-1}.
 \e*
 The details of approximation of derivatives with \reff{hermit} can be found in Lemma 2.1 in \cite{ftw}.

We intend to extend the result of \cite{ftw} to the nonlocal case. First observe that there is an obvious extension which could be done immediately by adding the following assumptions to Assumption {\bf F} in \cite{ftw}, i.e.

\no {\bf Assumption F}\quad {\it {\rm (i)} The nonlinearity $F$ is Lipschitz-continuous with respect to $(x,r,p,\gamma,\psi)$ uniformly in $t$, and $|F(\cdot,\cdot,0,0,0,0)|_\infty<\infty$;
\\
{\rm (ii)} $F$ is elliptic and dominated by the diffusion of the linear operator $\Lc^X$, i.e.
 \be\label{F2}
\nb_{\!\!\gamma} F\le a~&\mbox{on}&~\R^d\times\R\times\R^d\times\S_d\times\C_d;
 \ee
{\rm (iii)} $F_p\in{\rm Image}(F_\gamma)$ and $\big|F_p^\text{T}F_\gamma^{-}F_p\big|_\infty<+\infty$.
}

\no We remind that the  nonlocal nonlinearity $F$ is called elliptic if
\begin{enumerate}
\item $F$ is non-decreasing on the second derivative component, i.e.
\b*
F(t,x,r,p,\gamma_1,\psi) \le F(t,x,r,p,\gamma_2,\psi)
&~~\text{for}~~&
\gamma_1\le\gamma_2.
\e*
\item
$F$ is non-decreasing on the nonlocal  component, i.e.
\b*
F(t,x,r,p,\gamma,\psi_1) \le F(t,x,r,p,\gamma,\psi_2)
&~~\text{for}~~&
\psi_1\le\psi_2.
\e*
\end{enumerate}

\no Then we have the following Theorem.
\begin{Theorem}\label{thmobviousextension}
Let Assumption {\bf F} in hold true, and $|\mu|_1$, $|\sigma|_1<\infty$ and $\sigma$ is invertible. Also assume that the fully nonlinear PDE \reff{equationnl} has comparison for bounded functions. Then for every bounded Lipschitz function $g$, there exists a bounded function $v$ so that
 \b*
 v^h \longrightarrow v
 &&
 \mbox{locally uniformly}.
 \e*
In addition, $v$ is the unique bounded viscosity solution of problem \reff{equationnl}-\reff{terminalnl}.
\end{Theorem}
The proof is an straight forward implementation of the Subsection 3.2 of \cite{ftw}.

One of the major class of fully nonlinear PDE's, is the class of HJB equations which come from stochastic control problems arising in many application including finance. 
However, The nonlinearity of HJB equations do not satisfies Assumption {\bf F} in general. Even for local PDE's of HJB type, Assumption {\bf F} is not valid, because $F$ is not uniformly Lipschitz with respect to $x$. In addition, when the L\'evy measure $\nu$ is an infinite L\'evy measure, there is no chance for $F$ to be uniformly Lipschitz within respect to $\psi$. Therefore, we need to develop another theory for HJB equations.

The other problem which occurs in many applications is the lack of explicit form for nonlinearity $F$.
We present the following example in order to mention this problem.
\begin{Example}
Suppose we want to implement the scheme for the fully nonlinear equation of the form:
\b*
-v_t-F(x,Dv(t,x),D^2v(t,x),v(t,\cdot))&=&0
\\
v(T,\cdot)&=&g(\cdot),
\e*
where
 \be\label{example}
 F(x,p,\gamma,\psi)&:=&\sup_{\theta\in\R_+}\left \{\Lc^\theta(p,\gamma)+\int_{\R_*}\psi(x+\theta z)\nu(dz)\right \}\\
\Lc^\theta(p,\gamma)\label{diffop}
&:=&
\theta  b p+\frac12\theta^2a^2\gamma\\
\Ic(x,\psi)^\theta\label{intop}
&:=&
\int_{\R_*}\psi(x+\theta z)\nu(dz).
 \ee
This fully nonlinear equation solves the problem of portfolio management for one asset in the Black-Scholes model including jumps in asset price. For the sake of simplicity, for the moment we forget about infinite activity jumps.
Observe that if $\nu=0$ (the asset price do not jump) then $F$ becomes of the form:
 \b*
 F(x,p,\gamma,\psi)&:=&\sup_{\theta\in\R_+}\left\{\theta b p+\frac12\theta^2a^2\gamma\right \}.
 \e*
which could be given in explicit form by:
 \b*
 F(x,p,\gamma,\psi)&:=&-\frac{(b p)^2}{2 a^2\gamma},
 \e*
and the scheme could be easily implemented as in \cite{ftw} as well as more complicated examples. 
But, when $\nu\ne0$ (jump do exists), the explicit form for $F$ is not known and the supremum should be approximated. This problem is in common with other numerical methods for fully nonlinear PDE's  e.g. finite difference. Although his problem is obviously beyond the subject of this paper, we addressed it in this paper in order to mention that why we need to approximate the integral inside the supremum. 

The other problem, which appears in high dimensions, is the calculation of L\'evy integral inside supremum. Some numerical methods to approximate the supremum based on the calculation of the linear operator $\Lc^\theta+\Ic^\theta$ inside the supremum for different $\theta$s. Therefore, we proposed a Monte Carlo Quadrature method to approximate the integral in a purely probabilistic way. The MCQ could be considered independently in other applications.

From now on, we relax the assumption that $\nu$ is a finite measure. To be precise, we need to suppose that \reff{intop} is of the form
\b*
\Ic(x,\psi)
&:=&
\int_{\R_*}\left(\psi(x+\theta z)-\psi(x)-\1_{\{|z|\le1\}}\theta D\psi(x)\cdot z\right)\nu(dz).
\e*
In this case, there are two ways to work with singular L\'evy measure in numerical experiments; one is to truncate L\'evy measure near zero (as we also did for discretization of $X$) and the other is to approximate infinite small jumps by a Brownian motion. In both cases, the general form for the approximate $F$ is
 \b*
 F_{\kappa}(x,r,p,\gamma,\psi)&:=&\sup_{\theta\in\R_+}\left\{c_\kappa r+\theta b_\kappa p+\frac12\theta^2 a^2\gamma+\int_{\{|z|>\kappa\}}\psi(x+\theta z)\nu(dz)\right \}.
 \e*
where 
\b*
c_\kappa:=\int_{\{|z|>\kappa\}}\nu(dz)~~~\text{and}~~~b_\kappa:=b\int_{\{1\ge|z|>\kappa\}}z\nu(dz).
\e*
\end{Example}

We will introduce the modified scheme \reff{schemejump} in Section \ref{Sec:asymp} based on the approximation of nonlinearity $F$ obtained from truncation of infinite L\'evy measure and MCQ and then provide asymptotic results as in \cite{ftw} for nonlocal case.

The generalization of the result of \cite{ftw} for nonlocal PDE's would be easy if the function $F_\kappa$ were Lipschitz uniform on $\kappa$. But, for infinite L\'evy measures, this is never the case. To overcome this problem, we will show that $\kappa$ could be chosen dependent on $h$, so that the corresponding scheme satisfies the requirements of \cite{barlessouganidis} for the proof of convergence.

\section{Monte Carlo Quadrature ({\bf MCQ})}
\label{Sec:mcq}
In this section, we propose a Monte Carlo method the value of the following L\'evy generator:
\be\label{levyoperator}
\Ic[\varphi](x)
&:=&
\int_{\R_*^d}\left  (\varphi\left (x+\eta(z)\right )-\varphi(x)-\1_{\{|z|\le1\}}\eta(z)\cdot D\varphi(x)\right )\nu(dz).
\ee
 The method is  pure Monte Carlo method to approximate \reff{levyoperator} and, therefore could be used in the approximation of  L\'evy integral inside the scheme \reff{schemejump}.  Because, the result of this section is independent of the  numerical scheme \reff{schemejump} introduced in this paper, we organize this Section so that one can read it independently from other Section.\\
Through out this Section, we drop the dependency with respect to $(t,x)$ or other variables and for the sake of simplicity and just write $\eta(z)$. (For example in assumption {\bf F} in Section \ref{Sec:asymp} $\eta^{\alpha,\beta}(t,x,z)$ which depends on $(t,x,z,\alpha,\beta)$ will be considered as $\eta(z)$).\\
Notice that in order for \reff{levyoperator} to be well-defined for regular functions, we impose the following assumption on $\eta$:
\be
\label{eta}
\frac{|\eta(z)|}{|z|\wedge 1}
&\le&
 C,~~~\text{for some constant} ~~C.
\ee
We present MCQ in three cases with respect to the behavior of L\'evy measure near zero: 
\begin{itemize}
\item finite measure; $\int_{\{|z|\le1\}}\nu(dz)<\infty$,
\item infinite measure;
\begin{itemize}
\item case I: $\int_{\{|z|\le1\}}|\eta(z)|\nu(dz)<\infty$,
\item case II: $\int_{\{|z|\le1\}}|\eta(z)|^2\nu(dz)<\infty$.
\end{itemize}
\end{itemize}

\subsection{Finite L\'evy Measure}
When L\'evy measure is finite, we choose $\kappa=0$. In this case, we introduce Lemma \ref{lemintmc} which proposes a way to approximate the L\'evy integral of general form:
\be\label{integral}
\int_{\R_*^d}\varphi(x+\eta(z))\zeta(z)d\nu(z),
\ee
and then we use this Lemma to approximate the L\'evy infinitesimal generator \reff{levyoperator}.

Let $J$ be a jump Poisson measure with intensity given by L\'evy measure $\nu$, and $\{N_t\}_{t\ge0}$ be the Poisson process given by $N_t=\int_0^t\int_{\R_*^d}J(ds,dz) $ whose intensity is $\lambda:=\int_{\R_*^d}\nu(dz)$. By \reff{compoundpoisson}, we can write $\Xh^{x}$ by 
\be\label{compound1}
\Xh_t^{x}&=&x+\mu_0 t+\sigma W_t+\sum_{i=1}^{N_t}\eta(Z_i)
\ee
where $Z_i$s are i.i.d. random variables with law $\frac{1}{\lambda}\nu(dz)$.
We also introduce a  L\'evy process $Y_t$ by 
\be\label{compound2}
Y_t&=&\sum_{i=1}^{N_t}\zeta(Z_i).
\ee

Next Lemma shows that \reff{integral} could be approximated by a Monte Carlo formula purely free of integration.
\begin{Lemma}\label{lemintmc}
Let 
\be\label{approxintegral}
\hat\nu^{\eta,\zeta}_h(\varphi)(x):=\E\left[\int_{\R_*^d}\varphi(\Xh_h^x+\eta(z))\zeta(z)d\nu(z)\right].
\ee
Then,  for every bounded function $\varphi:\R^d\to\R$:
\b*
 \hat\nu^{\eta,\zeta}_h(\varphi)(x)
 &=&
 \frac{1}{h}\E[\varphi(\Xh^x_h)Y_h].
\e*
\end{Lemma}
\proof For the sake of simplicity, we just concentrate on the jump part of process $\Xh^{x}$ and without loss of generality, we write $\Xh_h^{x}=x+\sum_{i=1}^{N_h}\eta(Y_i)$. The right hand side can be expressed as:
\b*
\E\left[\varphi(\Xh^x_h)Y_h\right]
&=&
e^{-\lb h}\mathop{\sum}\limits_{n=0}^\infty \E\left[\varphi(\Xh^x_h)Y_h|N_h=n\right]\frac{(\lb h)^n}{n!}.
\e*
Then by \reff{compound1}-\reff{compound2},
\b*
\E\left[\varphi(\Xh^x_h)Y_h\right]
&=&
e^{-\lb h}\lb h\mathop{\sum}\limits_{n=1}^\infty \E\left[\varphi\left(x+\mathop{\sum}\limits_{i=1}^n \eta(Z_i)\right)\biggl(\mathop{\sum}\limits_{j=1}^n \zeta(Z_j)\biggr)\right]\frac{(\lb h)^{n-1}}{n!}\\
&=&e^{-\lb h}\lb h\mathop{\sum}\limits_{n=1}^\infty\frac{(\lb h)^{n-1}}{n!}\mathop{\sum}\limits_{j=1}^n \E\left[\varphi\left(x+\mathop{\sum}\limits_{i=1}^n \eta(Z_i)\right)\zeta(Z_j)\right].
\e*
Notice that in the above expression, the summation starts from $n=1$ because $Y_h=0$ when $N_h=0$. Because $Z_i$s are i.i.d. one can conclude that,
\b*
\mathop{\sum}\limits_{j=1}^n \E\left[\varphi\left(x+\mathop{\sum}\limits_{i=1}^n \eta(Z_i)\right)\zeta(Z_j)\right]
&=&
n\E\left[\varphi\left(x+\mathop{\sum}\limits_{i=1}^n \eta(Z_i)\right)\zeta(Z_1)\right]
\e*
Then, one can write
\b*
\E\left[\varphi\left(x+\eta(Z_1)\!\!+\!\!\mathop{\sum}\limits_{i=2}^n \eta(Z_i)\right)\zeta(Z_1)\right]
\!\!\!\!&\!\!=\!\!&\!\!\!\!
\E\left[\varphi\left(\eta(Z)+\hat X^x_h\right)\zeta(Z)|N_h=n-1\right],
\e*
where $Z$ is dependent of $Z_i$s but has the same law az $Z_i$s.
Therefore, we  can conclude that:
\b*
\E\left[\varphi(\Xh^x_h)Y_h\right]
&=&
e^{-\lb h}\lb h\mathop{\sum}\limits_{n=1}^\infty \E\left[\varphi(\eta(Z)+\hat X^x_h)\zeta(Z)|N_h=n-1\right]\frac{(\lb h)^{n-1}}{(n-1)!}.
\e*
But, we  know that
\b*
e^{-\lb h}\mathop{\sum}\limits_{n=1}^\infty \E\left[ \varphi(\eta(Z)+\hat X^x_h)\zeta(Z)|N_h=n-1\right]\frac{(\lb h)^{n-1}}{(n-1)!}
=
\E\left[\varphi(\eta(Z)+\hat X^x_h)\zeta(Z)\right]
\e*
Therefore,
\b*
\E\left[\varphi(\Xh^x_h)Y_h\right]
&=&
\lb h\E\left[\varphi(\eta(Z)+X^x_h)\zeta(Z)\right].
\e*
Because the density of $Z$ is $\frac{\nu(dz)}{\lb}$,
\b*
\E\left[\varphi( \Xh^x_h)Y_h\right]
&=&
h \E\left[\int_{\R^d_*}\varphi(\eta(z)+\Xh^x_h)\zeta(z)d\nu(z)\right].
\e*\ep\\
In the light of Lemma \reff{lemintmc}, we propose the following approximation for \reff{levyoperator}:
\b*
\Ic_h[\varphi](x):=\hat\nu_h^{\eta,1}-\varphi(x)\int_{\R^d_*}\nu(dz)- D\varphi(x)\cdot\int_{\R^d_*}\eta(z)\nu(dz).
\e*
Next Lemma provide error bound for this approximation.
 \begin{Lemma}\label{lemFapproxfinite}
 For any Lipschitz function $\varphi$ we have:
  \be\label{nonlinapprox}
 |(\Ic_{h}-\Ic)[\varphi]|_\infty
 &\le&
C\sqrt{h}|D\varphi|_\infty.
 \ee
 \end{Lemma}
 \proof
As a direct consequence of Lemma \reff{lemintmc}, $\hat\nu_h^{\eta,1}=\frac{1}{h}\E[\varphi(\Xh^x_h)N_h].$
Therefore, one can conclude that,
  \b*
  |(\Ic-\Ic_h)[\varphi]|_\infty
   &\le&
   C|D\varphi|_\infty\E\left [|\Xh_h^{x}-x|\right ].
  \e*
    So, because
\be\label{errorXh-x}
\E\left [|\Xh_h^x-x|\right ]
&\le&
C\left (h\int_{\R_*^d}|\eta(z)|\nu(dz)+\sqrt{h}\right ),
\ee    
which provides the result.\ep

\subsection{Infinite L\'evy Measure}
In the case of singular L\'evy measure, we  truncate L\'evy measure near zero and reduce the problem to a finite measure. In other words, for any $\kappa>0$ we have the truncation approximation of integral operator \reff{levyoperator}.
\b*
\Ic_\kappa[\varphi](x)&:=&\int_{\{|z|>\kappa\}} \!\!\!\!\!\!\!\!\!\!\!\!\!\!\left  (\varphi\left (x+\eta(z)\right )-\varphi(x)-\mathds{1}_{\{|z|\le1\}}\eta(z)\cdot D\varphi(x)\right )\nu(dz).
\e*
Then, we use  Lemma \reff{lemintmc} to present the MCQ approximation for \reff{levyoperator}.
 \b*
 \Ic_{\kappa,h}[\varphi](x):=\hat\nu^{\eta,1}_{\kappa,h}-\varphi(x)\!\!\!\int_{\{|z|>\kappa\}}\!\!\!\!\!\!\!\!\!\!\!\!d\nu(z)\!\!-\!\!\int_{\{1\ge|z|>\kappa\}}\!\!\!\!\!\!\!\!\!\!\!\eta(t,x,z)\cdot D\varphi(x)d\nu(z),
 \e*
where by Lemma \reff{lemintmc}
\b*
\hat\nu^{\eta,1}_{\kappa,h}:=\int_{\{|z|>\kappa\}}\varphi(\hat X^{x,\kappa}_h+\eta(t,x,z))\nu(dz)=h^{-1}\E\left[\varphi(\hat X^{x,\kappa}_h)N_h^\kappa\right]
\e*
Following Lemma provides the error of MCQ approximation of \reff{levyoperator} in the case of infinite L\'evy measure.

\begin{Lemma}
\label{lemFapproxinfinite}
Let function $\varphi$ be Lipschitz.
\begin{enumerate}
\item If $\int_{\{|z|\le1\}}|z|\nu(dz)<\infty$, then
\be\label{nonlinapprox1}
 |(\Ic_{\kappa,h}-\Ic)[\varphi]|_\infty
 \!\!&\!\!\le\!\!&\!\!C|D\varphi|_\infty\left (\sqrt{h}+\int_{\{0<|z|\le\kappa\}}\!\!\!\!\!\!\!\!\!\!\!\!|z|\nu(dz)\right ).
 \ee
\item If $\int_{\{|z|\le1\}}|z|^2\nu(dz)<\infty$, then
\be\label{nonlinapprox2}
 |(\Ic_{\kappa,h}-\Ic)[\varphi]|_\infty
 \!\!&\!\!\le\!\!&\!\!C\biggl(|D\varphi|_\infty\Bigr(\sqrt{h}+h\!\!\int_{\{|z|>\kappa\}}\!\!\!\!\!\!\!\!\!\!\!\!|z|\nu(dz)\Bigr)+|D^2\varphi|_\infty\!\!\int_{\{0<|z|\le\kappa\}}\!\!\!\!\!\!\!\!\!\!\!\!\!\!\!|z|^2\nu(dz)\biggr).
 \ee
\end{enumerate}
 \end{Lemma}
 \proof
\begin{enumerate}
\item Notice that,
\b*
|(\Ic-\Ic_{\kappa,h})[\varphi]|_\infty\le|(\Ic-\Ic_{\kappa})[\varphi]|_\infty+|(\Ic_{\kappa}-\Ic_{\kappa,h})[\varphi]|_\infty.
\e*
By \reff{eta}, the truncation error is given by:
\be\label{truncerror1}
|(\Ic-\Ic_{\kappa})[\varphi]|_\infty\le2|D\varphi|_\infty\int_{\{0<|z|\le\kappa\}}|\eta(z)|\nu(dz).
\ee
On the other hand, by \reff{errorXh-x} and \reff{eta}, we observe that
\b*
|(\Ic_{\kappa}-\Ic_{\kappa,h})[\varphi]|_\infty
&\le& 
C|D\varphi|_\infty\left (h\int_{\{|z|>\kappa\}}|\eta(z)|\nu(dz)+\sqrt{h}\right )\\
&\le&
C|D\varphi|_\infty\left (h\int_{\{|z|>\kappa\}}|z|\nu(dz)+\sqrt{h}\right )
\e*
which together with \reff{truncerror1} provides the result.
\item
By \reff{eta}, the truncation error is given by:
\be\label{truncerror2}
|(\Ic-\Ic_{\kappa})[\varphi]|_\infty\le C|D^2\varphi|_\infty\int_{\{0<|z|\le\kappa\}}|z|^2\nu(dz),
\ee
for any function $\varphi$ with bounded derivatives up to second order. 
On the other hand, \reff{errorXh-x} allows us to calculate the Monte Carlo error by:
\b*
|(\Ic_{\kappa}-\Ic_{\kappa,h})[\varphi]|_\infty
&\le& 
C|D\varphi|_\infty\left (h\int_{\{|z|>\kappa\}}|z|\nu(dz)+\sqrt{h}\right )
\e*
which completes the proof.\ep
\end{enumerate}

\section{Asymptotic results}
\label{Sec:asymp}
This section is devoted to the convergence result for the scheme \reff{schemejump}. 
We first remind the notion of viscosity solution and provide the assumptions required for the main results together with the statement of main results. Then we provide the proof of the results in two following subsection.

we need to impose the following assumption on the nonlinearity $F$ to obtain the convergence Theorem.
\vspace*{5mm}

\noindent{\bf Assumption IHJB1}: {\it Function $F$ satisfies:
\b*
 \frac12a(t,x)\cdot\gamma+\mu(t,x)\cdot p+F(t,x,r,p,\gamma,\psi)
\!\!\! &\!\!:=\!\!&\!\!\!
 \inf_{\alpha\in\Ac}\sup_{\beta\in\Bc}\left\{\Lc^{\alpha,\beta}(t,x,r,p,\gamma)+\Ic^{\alpha,\beta}(t,x,r,p,\gamma,\psi)\right\}
 \e*
for given sets $\Ac$ and $\Bc$ where
   \b*
 \Lc^{\alpha,\beta}(t,x,r,p,\gamma)
\!\!\! &\!\!:=\!\!&\!\!\!
 \frac{1}{2}
 a^{\alpha,\beta}(t,x)\cdot\gamma
 +b^{\alpha,\beta}(t,x)\cdot p+c^{\alpha,\beta}(t,x)r+k^{\alpha,\beta}(t,x),
 \e*
and
 \b*
 \Ic^{\alpha,\beta}(t,x,r,p,\psi)
\!\!\! &\!\!:=\!\!&\!\!\!
\int_{\R^d_*}\left(\psi\left(x+\eta^{\alpha,\beta}(t,x,z)\right)-r-\mathds{1}_{\{|z|\le1\}}\eta^{\alpha,\beta}(t,x,z)\cdot p\right)\nu(dz)
 \e*
where for any $(\alpha,\beta)\in\Ac\times\Bc$, $a^{\alpha,\beta}$, $b^{\alpha,\beta}$, $c^{\alpha,\beta}$, $k^{\alpha,\beta}$ and $\eta^{\alpha,\beta}$  satisfy
\b*
\sup_{\alpha\in\Ac, \beta\in\Bc}\biggl\{|a^{\alpha,\beta}|_1+|b^{\alpha,\beta}|_1+|c^{\alpha,\beta}|_1+|k^{\alpha,\beta}|_1+
\frac{|\eta^{\alpha,\beta}(\cdot,z)|_1}{|z|\wedge1}\biggr\}<\infty.
\e*
The nonlinearity is dominated by the diffusion of the linear operator $\Lc^X$, i.e. for any $t$, $x$, $z$, $\alpha$ and $\beta$
 \begin{align}
\label{F1} |a^{-}\cdot a^{\alpha^*,\beta^*}|_1<\infty~~~\text{and}~~~&0 \le a^{\alpha,\beta}\le a, \\
\label{F2}\eta^{\alpha,\beta},~b^{\alpha,\beta}\in{\rm Image}(a^{\alpha,\beta})~~~\text{and}~~~&\sup_{\alpha\in\Ac, \beta\in\Bc}|(b^{\alpha,\beta})^{\rm T}(a^{\alpha,\beta})^{-} b^{\alpha,\beta}|_\infty<\infty,\\
\nonumber&\sup_{\alpha\in\Ac, \beta\in\Bc}\frac{|(\eta^{\alpha,\beta})^{\rm T}(a^{\alpha,\beta})^{-} b^{\alpha,\beta}|_\infty}{1\wedge |z|}<\infty\\
\nonumber&\sup_{\alpha\in\Ac, \beta\in\Bc}\frac{|(\eta^{\alpha,\beta})^{\rm T}(a^{\alpha,\beta})^{-} \eta^{\alpha,\beta}|_\infty}{1\wedge |z|^2}<\infty.
 \end{align}
}
\begin{Remark}
 {\rm
A function $F$ which satisfies Assumption {\bf IHJB1} is not well-defined for arbitrary $(t,x,r,p,\gamma,\psi)$ $\in\R_+\times\R^d\times\R\times\R^d\times\S_d\times\Cc_d$. But, for any second order differentiable function,  $\psi$,  with bounded derivatives with respect to $x$, $F(t,x,\psi(t,x),D\psi(t,x),D^2\psi(t,x),\psi(t,\cdot))$ is well-defined. 
}
\end{Remark}

Now, we propose a Monte Carlo scheme for \reff{equationnl}-\reff{terminalnl} based on the same idea as in \cite{ftw}, and also the approximation of the nonlinearity.
 \be\label{schemejump}
 v^{\kappa,h}(T,.)=g
 &\mbox{and}&
 v^{\kappa,h}(t_i,x)=\TT_{\kappa,h}[v^{\kappa,h}](t_i,x),
 \ee
where for every function $\psi:\R_+\times\R^d\longrightarrow\R$ with exponential growth:
 \be\label{TTkappah}
 \TT_{\kappa,h}[\psi](t,x)\!\!\!\!\!&:=&\!\!\!\!\!\E\left[\psi\left(t+h,\hat X_h^{t,x,\kappa}\right)\right]
                   +hF_{\kappa,h}\left(t,x,\Dc_h\psi,\psi(t+h,\cdot)\right),\\
\nonumber&&~\Dc_h\psi:=\left(\Dc_h^0\psi,\Dc_h^1\psi,\Dc_h^2\psi\right),
 \ee
\b*
{F}_{\kappa,h}(t,x,r,p,\gamma,\psi)
\!\!\!&\!\!=\!\!&\!\!\!\inf_{\alpha\in\Ac}\sup_{\beta\in\Bc}\biggl\{ 
\frac{1}{2}
 a^{\alpha,\beta}(t,x)\cdot\gamma+b^{\alpha,\beta}(t,x)\cdot p+c^{\alpha,\beta}(t,x)r+k^{\alpha,\beta}(t,x)\\
\!\!\!&\!\!\!\!&\!\!\! \hspace*{1.5cm}+\int_{\{|z|\ge\kappa\}}\!\!\biggl(\hat\nu^{\eta^{\alpha,\beta},1}_h(\psi(t,\cdot))(x)-r-\eta^{\alpha,\beta}(t,x,z)\cdot p\biggr)\nu(dz)
\biggr\},
\e*
and
 \be\label{hermit}
 \Dc_h^k\psi(t,x)
 &:=&
 \E\left[\psi(t+h,\hat X^{t,x,\kappa}_h)H^h_k(t,x)\right], k=0,1,2,
 \ee
 where
  \b*
 H^h_0=1,
 &H^h_1=\left({\sigma^{\rm T}}\right)^{-1}\;\Frac{W_h}{h},&
 H^h_2=\left({\sigma^{\rm T}}\right)^{-1}\;\frac{W_hW^{\rm T}_h-h\mathbf{I}_d}{h^2}\;\sigma^{-1}.
 \e*
 The details of approximation of derivatives with \reff{hermit} can be found in Lemma 2.1 in \cite{ftw}. 
 In order to have the convergence result, we also need to impose the following assumption over $F_{\kappa,h}$.

\noindent{\bf Assumption Inf-Sup}: {\it For any $\kappa>0$, $t\in[0,T]$, $x$ and $x'\in\R^d$  and any Lipschitz functions $\psi$ and $\varphi$, there exists a $(\alpha^*,\beta^*)\in\Ac\times\Bc$ such that 
\b*
\Phi_\kappa^{\alpha^*,\beta^*}[\psi,\varphi](t,x,x')=\Jc_\kappa^{\alpha^*,\beta^*}[\psi](t,x)- \Jc_\kappa^{\alpha^*,\beta^*}[\varphi](t,x')
\e*
where
\be\label{Phi}
\Phi_\kappa^{\alpha,\beta}[\psi,\varphi](t,x):=\inf_\alpha \Jc_\kappa^{\alpha,\beta}[\psi](t,x)-\sup_\beta \Jc_\kappa^{\alpha,\beta}[\varphi](t,x'),
\ee
and
\b*
\Jc_\kappa^{\alpha,\beta}[\phi](t,x)
&:=&
\frac{1}{2}
 a^{\alpha,\beta}\cdot D^2\phi(t,x)+b^{\alpha,\beta}\cdot D\phi(t,x)+c^{\alpha,\beta}\phi(t,x)+k^{\alpha,\beta}(t,x)\\
\!\!\!&\!\!\!\!&\!\!\! \hspace*{1.5cm}+\int_{\{|z|\ge\kappa\}}\!\!\!\!\!\!\!\bigl(\hat\nu^{\eta^{\alpha,\beta},1}_h(\phi(t,\cdot))(x)-\phi(t,x)-\eta^{\alpha,\beta}(t,x,z)\cdot D\phi(t,x)\bigr)\nu(dz).
\e*
}
%
%

The first result concerns the convergence of the convergence of $v^{\kappa,h}$ for $\kappa$ appropriately chosen with respect to $h$.

\begin{Theorem}[Convergence]\label{thmconvjump}
Let $\eta$, $\mu$ and $\sigma$ be bounded and Lipschitz continuous on $x$ uniformly on $t$ and $z$, $\sigma$ is invertible and Assumptions {\bf IHJB1} and {\bf Inf-Sup} hold true, and assume that \reff{equationnl} has comparison for bounded functions. Then, if  $\kappa_h$ is such that:
\be\label{thetakappah}
\lim_{h\to0}\kappa_h=0~~~&\text{and}&~~~\limsup_{h\to0}\theta_{\kappa_h}^2 h=0
\ee
where
\be\label{thetakappa}
\theta_\kappa&:=&\sup_{\alpha,\beta}|\theta_\kappa^{\alpha,\beta}|_\infty,
\ee
with
\b*
\theta_\kappa^{\alpha,\beta}&:=&c^{\alpha,\beta}+\int_{\{|z|\ge\kappa\}}\!\!\!\!\!\!\!\!\!\!\!\!\nu(dz)+\frac{1}{4}\biggl(b^{\alpha,\beta}\!\!-\!\!\int_{\{1>|z|\ge\kappa\}}\!\!\!\!\!\!\eta^{\alpha,\beta}(z)\nu(dz)\biggr)^\text{T}\\
\nonumber&&\times({a^{\alpha,\beta}})^{-}\biggl(b^{\alpha,\beta}\!\!-\!\!\int_{\{1>|z|\ge\kappa\}}\!\!\!\!\!\!\eta^{\alpha,\beta}(z)\nu(dz)\biggr),
\e*
then $v^{\kappa_h,h}$ converges to some function $v$ locally uniform. In addition, $v$ is the unique viscosity solution of \reff{equationnl}-\reff{terminalnl}.\\
Specially, if L\'evy measure is finite for the choice of $\kappa_h=0$ the assertion of the Theorem hold true.
\end{Theorem}

\begin{Remark}\label{hkappah}
{\rm
It is always possible to choose $\kappa_h$ such that \reff{thetakappah} is satisfied. To see this, notice that $\theta_\kappa$ in \reff{thetakappa} is  non-increasing on $\kappa$
\b*
\lim_{\kappa\to0}\theta_\kappa=+\infty~~\text{and}~~\limsup_{\kappa\to\infty}\theta_\kappa<\infty.
\e*
Then, we define $\kappa_h:=\inf\{\kappa|\theta_\kappa\le h^{-\frac12}\}+h$. By the definition of $\kappa_h$, $\theta_{\kappa_h}\le h^{-\frac12}$. Because
Observe that $\kappa_h$ is non-decreasing with respect to $h$ and $\lim_{h\to0}\kappa_h=0$. \\
If there exists a $q$ such that, $q:=\lim_{h\to0}\kappa_h>0$, then, for $\kappa<q$, we would have $\theta_{\kappa}=\infty$ which obviously contradicts the fact that for $\kappa>0$, $\theta_\kappa<\infty$. Therefore, $\kappa_h$ satisfies \reff{thetakappa}.
}
\end{Remark}

\begin{Remark}{\rm
The choice of $\kappa_h$ in the above Theorem seems to be crucial for the convergence. Otherwise, we only have the following convergence result.}
\end{Remark}
\begin{Proposition}
Under the same assumption as Theorem \ref{thmconvjump},
when L\'evy measure $\nu$ is infinite, for every Lipschitz bounded function $g$, we have
\b*
\lim_{\kappa\to0}\lim_{h\to0}v^{\kappa,h}=v
\e*
where $v$ is the unique viscosity solution of \reff{equationnl}-\reff{terminalnl} assuming that it exists.
\end{Proposition}
\proof.
Let $v^\kappa$ be the solution of the following problem:
\be
 \label{equationnlkappa}
\!\!\!\!\!\!\!\!\!\!\!\! &&\!\!\!\!\!\!\!\!\!\!\!\!-\Lc^X \!v^\kappa(t,x) \!\!-\!\! F_\kappa\!\!\left(t,x,\!v^\kappa\!(t,x),\!D v^\kappa\!(t,x),\!D^2v^\kappa\!(t,x),\!v^\kappa\!(t,\cdot)\right)= 0, \mbox{on}~\!\![0,T)\!\!\times\!\!\R^d,\\
\!\!\!\!\!\!\!\!\!\!\!\!&&\!\!\!\!\!\!\!\!\!\!\!\!v^\kappa(T,\cdot)=g(\cdot), ~~~~~~~~~~~~~~~~~~~~~~~~~~~~~~~~~~~~~~~~~~~~~~~~~~~~~~~~~~~~~~~~\mbox{on}~\in\R^d. \label{terminalnlkappa}
\ee
where 
$F_\kappa:\R_+\times\R^d\times\R\times\R^d\times\Sc_d\times\Cc_d\rightarrow\R$ is given by:
\b*
 F_\kappa(t,x,r,p,\gamma,\psi)
 &:=&
 \inf_{\alpha\in\Ac}\sup_{\beta\in\Bc}\left\{\Lc^{\alpha,\beta}(t,x,r,p,\gamma)+\Ic_\kappa^{\alpha,\beta}(t,x,r,p,\gamma,\psi)\right\}
 \e*
where
 \be\label{vkappav}
 \Ic_\kappa^{\alpha,\beta}(t,x,r,p,\gamma,\psi)\!\!&\!\!\!\!:=\!\!\!\!&\!\!\int_{\{|z|\ge\kappa\}}\!\!\!\!\!\!\!\!\!\!\!\!\!\bigl(\psi\bigl(x+\eta^{\alpha,\beta}(t,x,z)\bigr)-r-\mathds{1}_{\{|z|\le1\}}\eta^{\alpha,\beta}(t,x,z)\cdot p\bigr)\nu(dz)~~~~~~~~~~~~~
 \ee
where $a^{\alpha,\beta}$, $b^{\alpha,\beta}$, $c^{\alpha,\beta}$, $k^{\alpha,\beta}$ and $\eta^{\alpha,\beta}$ are as in Assumption {\bf IHJB1}.
Let $v^{\kappa,h}$ be the approximate solution given by the scheme \reff{schemejump}. Let $\kappa>0$ be fixed. Because the truncated L\'evy measure is finite, by Theorem \ref{thmconvjump}, $v^{\kappa,h}$ converges to $v^\kappa$ locally uniformly as $h\to0$. 
Let $v^\kappa$ be the solution of \reff{equationnlkappa}-\reff{terminalnlkappa}.
By Theorem 5.1 of \cite{biswasjakobsenkarlsen2} and Assumption {\bf F}, we have:
\be\label{vkappav}
|v-v^\kappa|_\infty&\le&C\sup_{\alpha,\beta}\left\{\left (\int_{0<|z|<\kappa}| \eta^{\alpha,\beta}(\cdot,z)|^2_\infty\nu(dz)\right )^\frac12\right\}\\
&\le&C\left (\int_{0<|z|<\kappa}|z|^2_\infty\nu(dz)\right )^\frac12.
\ee
Therefore, one can choose $\kappa>0$ so that $|v^\kappa-v|_\infty$ be small enough. Then, when $h$ goes to $0$, $v^{\kappa,h}$ converges to $v^\kappa$.\ep\\

\no The above limit proposes to implement the numerical scheme in two steps:
\begin{itemize}
\item First by choosing $\kappa$ so that $v^{\kappa}$ is near enough to $v$, we obtain a uniform approximation of $v$.
\item Second by sending $h\to0$, we obtain locally uniform convergence of $v^{\kappa,h}$ to $v^\kappa$.
\end{itemize}
Notice that the above convergence is not uniformly on $(\kappa,h)$. However, the convergence in Theorem \ref{thmconvjump}, is uniform on $h$ when the choice of $\kappa$ is made suitably dependent on $h$.

\vspace*{5mm}

\begin{Remark}\label{generalfinal}{\rm
By Remark 3.7 in \cite{ftw},the boundedness condition on $g$ can be relaxed.}
\end{Remark}
In order to obtain the rate of convergence result, we impose Assumptions {\bf IHJB2} and {\bf IHJB2+} which restrict us to concave nonlinearities.\\

\no {\bf Assumption IHJB2}\quad {\it The nonlinearity $F$  satisfies Assumption {\bf IHJB1} with $\Bc$ be a singlton set
}\\
\begin{Remark}\label{concavenonlinearity}{\rm
Therefore, when the nonlinearity $F$ satisfies {\bf IHJB}, we can drop the super script $\beta$ and write $F$  by
 \b*
 F(t,x,r,p,\gamma,\psi)
 &:=&
 \inf_{\alpha\in\Ac}\left\{\Lc^{\alpha}(t,x,r,p,\gamma)+\Ic^\alpha(t,x,r,p,\gamma,\psi)\right\}
 \e*
  where
  \b*
 \Lc^{\alpha}(t,x,r,p,\gamma)
 &:=&
 \frac{1}{2}
 {\Tr{(a^\alpha)^{\rm T}}(t,x)\gamma}
 +b^{\alpha}(t,x)p+c^{\alpha}(t,x)r+k^\alpha(t,x),
 \e*
and
 \b*
 \Ic^\alpha(t,x,r,p,\psi)&:=&\int_{\R^d_*}\left(\psi\left(x+\eta^\alpha(t,x,z)\right)-r-\mathds{1}_{\{|z|\le1\}}\eta^\alpha(t,x,z)\cdot p\right)\nu(dz).
 \e*
In this case, the nonlinearity is a concave function of $(r,p,\gamma,\psi)$.
}
\end{Remark}
\no {\bf Assumption IHJB+}\quad {\it The nonlinearity $F$ satisfies {\bf IHJB2} and for any $\delta>0$, there exists a finite set $\left\{\alpha_i\right\}_{i=1}^{M_\delta}$ such that for any $\alpha\in\Ac$:
\b*
\mathop{\inf}\limits_{1\le i\le M_\delta}\!\biggl\{|\sigma^\alpha-\sigma^{\alpha_i}|_\infty\!+\!|b^\alpha-b^{\alpha_i}|_\infty\!+\!|c^\alpha-c^{\alpha_i}|_\infty\!+\!|k^\alpha-k^{\alpha_i}|_\infty\!+\!\!\!\mathop{\int}\limits_{\R^d_*}\!\!|(\eta^\alpha-\eta^{\alpha_i})(\cdot,z)|_\infty^2d\nu(z)\!\biggr\}\!\le\!\delta.
\e*
}

\begin{Remark}{\rm
The Assumption IHJB+ is satisfied if $\Ac$ is a compact separable topological space and $\sigma^\alpha(\cdot)$, $b^\alpha(\cdot)$, and $c^\alpha(\cdot)$   are continuous maps from $\Ac$ to $C_b^{\frac{1}{2},1}([0,T]\times\R^d)$; the space of bounded maps which are Lipschitz on $x$ and $\frac{1}{2}$-H\"older on $t$ and $\eta^\alpha(\cdot)$ is continuous maps from $\Ac$ to $\Bigl\{\varphi:[0,T]\times\R^d\times\R^d_*\to\R\Bigl|\int_{R^d_*}|\varphi(\cdot,z)|_\infty^2\nu(dz)<\infty\Bigr\}$.
}
\end{Remark}

\begin{Theorem}[Rate of Convergence]\label{thmrateconvjump}
Assume that the final condition $g$ is bounded and Lipschitz-continuous.
Then, there is a constant $C>0$ such that 
\begin{itemize}
\item under Assumption IHJB,\\
 $v-v^{\kappa,h} \le C\left (h^\frac{1}{4}+h\theta_{\kappa}^2+h\eps^{-3}+h^{\frac{3}{4}}\theta_\kappa+h\sqrt{\theta_\kappa}+h^{-\frac{1}{4}}\int_{\{|z|\le\kappa\}}\!\!|z|^2\nu(dz)\right )$.
\item under Assumption IHJB+, \\
$-C\left (h^{1/10}+h^{\frac{7}{10}}\theta_\kappa+h\sqrt{\theta_\kappa}+h^{-\frac{3}{10}}\int_{\{|z|\le\kappa\}}\!\!|z|^2\nu(dz)\right )\le v-v^{\kappa,h}$.
\end{itemize}
In addition, if it is possible to find $\kappa_h$ such that
\be \label{kappahrate}
\lim_{h\to0}\kappa_h=0,~~\limsup_{h\to0}h^\frac{3}{4}\theta_{\kappa_h}^2<\infty~~\text{and}~~
\limsup_{h\to0}h^{-\frac{1}{2}}\int_{0<|z|<\kappa_h}\!\!|z|^2\nu(dz)<\infty,
\ee
then, there is a constant $C>0$ such that 
\begin{itemize}
\item under Assumption IHJB, $v-v^{\kappa_h,h} \le C h^{1/4}$.
\item under Assumption IHJB+, $-C h^{1/10}\le v-v^{\kappa_h,h}$.
\end{itemize}
\end{Theorem}
Next example shows the case where the conditions of the above Theorem on the choice of $\kappa$ is satisfied. If it is not the case in some situations, it does mean that the rate of convergence is less than what is proposed by Theorem \ref{thmrateconvjump}.
\begin{Example}\label{kappahchoice}
For the L\'evy measue 
\b*
\nu(dz)=\1_{\R^d_*}|z|^{-d-1}dz,
\e*
one can always find $\kappa_h$ such that the condition of Theorem \ref{thmrateconvjump} is satisfied. In the other words, it is always enough to choose $\kappa_h$ such that
\b*
\limsup_{h\to0}h^{-\frac{1}{2}}\kappa_h<\infty.
\e*
\end{Example}

\subsection{Convergence}
We suppose the all the assumptions of Theorem \ref{thmconvjump} holds true throughout this subsection.

We first manipulate the scheme to provide strict monotonicity by the similar idea as in Remark 3.13 and Lemma 3.19 in \cite{ftw}. Let $u^{\kappa,h}$ be the solution of 
\be\label{overlinescheme}
 u^{\kappa,h}(T,\cdot)=g
 &\text{and}&
 u^{\kappa,h}(t_i,x)=\overline{\TT}_{\kappa,h}[u^{\kappa,h}](t_i,x),
 \ee
where
 \be\label{overlineTT}
 \overline{\TT}_{\kappa,h}[\psi](t,x)\!\!\!\!\!&:=&\!\!\!\!\!\E\left[\psi\left(t+h,\hat X_h^{t,x,\kappa}\right)\right]
                   +h\overline{F}_{\kappa,h}\left(t,x,\Dc_h\psi,\psi(t+h,\cdot)\right)
 \ee
 and
 \b*
\overline{F}_{\kappa,h}(t,x,r,p,\gamma,\psi)
\!\!\!&\!\!=\!\!&\!\!\!\sup_{\alpha}\inf_{\beta}\biggl\{ 
\frac{1}{2}
 a^{\alpha,\beta}\cdot\gamma+b^{\alpha,\beta}\cdot p+(c^{\alpha,\beta}+\theta_\kappa)r+e^{\theta_\kappa(T-t)}k^{\alpha,\beta}(t,x)\\
\!\!\!&\!\!\!\!&\!\!\! +\int_{\{|z|\ge\kappa\}}\!\!\!\!\!\!\!\bigl(\hat\nu_{\kappa,h}^{\eta^{\alpha,\beta},1}(\psi)-r-\mathds{1}_{\{|z|\le1\}}\eta^{\alpha,\beta}(z)\cdot p\bigr)\nu(dz)
\biggr\}.
\e*

\begin{Remark}\label{reminfsupbar}
 {\rm 
Assumption {\bf Inf-Sup} is also true if we replace $\Jc^{\alpha,\beta}_\kappa$ by
\b*
\bar\Jc^{\alpha,\beta}_\kappa[\psi](t,x)
\!\!\!&\!\!\:=\!\!&\!\!\!
\frac{1}{2}
 a^{\alpha,\beta}\cdot D^2\phi(t,x)+b^{\alpha,\beta}\cdot D\phi(t,x)+(c^{\alpha,\beta}+\theta_\kappa)\phi(t,x)+e^{\theta_\kappa(T-t)}k^{\alpha,\beta}(t,x)\\
&& \hspace*{1cm}+\int_{\{|z|\ge\kappa\}}\!\!\!\!\!\!\!\bigl(\hat\nu^{\eta^{\alpha,\beta},1}_h(\phi(t,\cdot))(x)-\phi(t,x)-\eta^{\alpha,\beta}(t,x,z)\cdot D\phi(t,x)\bigr)\nu(dz).
\e*
The proof is straight forward.
}
\end{Remark}

We have the following Lemma which shows that for proper choice of $\theta_\kappa$ the scheme \reff{overlinescheme} is strictly monotone.

\begin{Lemma}\label{lemmonjumpkappa}
Let $\theta_\kappa$ be as in  \reff{thetakappa}
and $\varphi$ and $\psi:~[0,T]\times\R^d\longrightarrow\R$ be two bounded functions. Then:
 \b*
 \varphi\le\psi
 &\Longrightarrow&
 \overline{\TT}_{\kappa,h}[\varphi]\le \overline{\TT}_{\kappa,h}[\psi].
 \e*
\end{Lemma}
\proof Let $f:=\psi-\varphi\ge 0$ where $\varphi$ and $\psi$ are as in the statement of the lemma. For simplicity, we drop the dependence on $(t,x)$ when it is not necessary. By Assumption {\bf IHJB1} and Lemma \reff{lemintmc}, we can write:
\b*
\overline{\TT}_{\kappa,h}[\psi]-\overline{\TT}_{\kappa,h}[\varphi]
\!\!&\!\!=\!\!&\!\!
\E[f(t+h,\Xh_h)]\\
\!\!&\!\!\!\!&\!\!+h\left(\inf_{\alpha}\sup_{\beta}
\bar\Jc_\kappa^{\alpha,\beta}[\what\psi](t+h,x)-\inf_{\alpha}\sup_{\beta}\bar\Jc_\kappa^{\alpha,\beta}[\what\varphi](t+h,x)\right),
\e*
where $\what\phi(t,x):=\E[\phi(t,\Xh_h^x)]$ for $\phi=\varphi$ or $\psi$.
Therefore,
\b*
\overline{\TT}_{\kappa,h}[\psi]-\overline{\TT}_{\kappa,h}[\varphi]
&\ge&\E[f(t+h,\Xh_h)]+h\bar\Phi_\kappa^{\alpha,\beta}[\what\psi,\what\varphi](t+h,x,x),
\e*
where $\bar\Phi_\kappa^{\alpha,\beta}$ is defined by
\b*
\bar\Phi_\kappa^{\alpha,\beta}[\psi,\varphi](t,x):=\inf_\alpha \bar\Jc_\kappa^{\alpha,\beta}[\psi](t,x)-\sup_\beta \bar\Jc_\kappa^{\alpha,\beta}[\varphi](t,x').
\e*
By Asumption {\bf Inf-Sup}, there exists $(\alpha^*,\beta^*)$ so that
\b*
\overline{\TT}_{\kappa,h}[\psi]-\overline{\TT}_{\kappa,h}[\varphi]
&\ge&\E[f(t+h,\Xh_h)]+h\left(\bar\Jc_\kappa^*[\what\psi](t+h,x)-\bar\Jc_\kappa^*[\what\varphi](t+h,x)\right).
\e*
Observe that by the linearity of $\bar\Jc_\kappa^{\alpha,\beta}$, one can write:
\b*
\bar\Jc_\kappa^{\alpha,\beta}[\hat \phi](t+h,x)
&=&
\E\left[
\bar\Jc_\kappa^{\alpha,\beta}[\phi](t+h,\Xh_h)\right].
\e*
By the definition of $\bar\Jc_\kappa^{\alpha,\beta}$ and Lemma 2.1 in \cite{ftw},
\begin{align*}
 \overline{\TT}_{\kappa,h}[\psi]-&\overline{\TT}_{\kappa,h}[\varphi]
\ge
\E\biggl[f(\hat X_h)\biggl(1+h\Bigl(c^{\alpha^*,\beta^*}_\kappa+\theta_\kappa
+b^{\alpha^*,\beta^*}_\kappa\cdot{(\sigma^{\rm T}})^{-1}\frac{W_h}{h}\\
& + \frac{1}{2} a^{\alpha^*,\beta^*}\cdot {(\sigma^{\rm T}})^{-1}\frac{W_h W^\text{T}_h-h\text{I}_{d}}{h^2}\sigma^{-1}\Bigr)\biggr)+h\hat\nu_h^{\eta^{\alpha^*,\beta^*},1}(f)\biggr],
\end{align*}
where  $b_\kappa^{\alpha,\beta}=b^{\alpha,\beta}-\int_{\{1>|z|\ge\kappa\}}\!\!\eta^{\alpha,\beta}(z)\nu(dz)$ and $c_\kappa^{\alpha,\beta}=c^{\alpha,\beta}-\int_{\{|z|\ge\kappa\}}\!\!\nu(dz)$.\\

\no Therefore, by the same argument as in Lemma 3.12 in \cite{ftw}, one can write:
\begin{align*}
\overline{\TT}_{\kappa,h}[\psi]-\overline{\TT}_{\kappa,h}[\varphi]
&\ge
\E\biggl[f(\hat X_h)\biggl(1-\frac12 a^{\alpha^*,\beta^*}\cdot a^{-1} +h\Bigl(|A_h^{\alpha^*,\beta^*}|^2+c_\kappa^*+\theta_\kappa\\
&-\frac{1}{4} ({b^{\alpha^*,\beta^*}_\kappa})^\text{T}({a^{\alpha^*,\beta^*}})^{-}b^{\alpha^*,\beta^*}_\kappa\Bigr)\biggr) +h\hat\nu_h^{\eta^{\alpha^*,\beta^*},1}(f)\biggr],
\end{align*}
where
 \be\label{Ah}
 A^{\alpha^*,\beta^*}_h &:=& \frac{1}{h}{(\sigma^{\alpha^*,\beta^*}})^{1/2}({\sigma^{\rm T}})^{-1}W_h
                +\frac{1}{2}(({\sigma^{\alpha^*,\beta^*}})^{-})^{1/2}b^{\alpha^*,\beta^*}_\kappa.
 \ee
Therefore, by positivity of $f$ and Assumption {\bf IHJB1}, one can deduce:
\b*
\overline{\TT}_{\kappa,h}[\psi]-\overline{\TT}_{\kappa,h}[\varphi]
&\ge&
h\E\biggl[f(\hat X_h)\biggl(c_\kappa^*+\theta_\kappa-\frac{1}{4} ({b^{\alpha^*,\beta^*}_\kappa})^\text{T}({a^{\alpha^*,\beta^*}})^{-}b^{\alpha^*,\beta^*}_\kappa\biggr)\biggr]
\e*
By the choice of $\theta_\kappa$ in \reff{thetakappa},  we have
\b*
\overline{\TT}_{\kappa,h}[\psi]-\overline{\TT}_{\kappa,h}[\varphi]
&\ge&
0.
\e*
Then, sending $\eps$ to zero provides the result.\ep
\vspace*{5mm}

The following Corollary shows the monotonicity of scheme \ref{schemejump}. 
\begin{Corollary}\label{cormonjump}
Let $\varphi, \psi:~[0,T]\times\R^d\longrightarrow\R$ be two bounded functions. Then:
 \b*
 \varphi\le\psi
 &\Longrightarrow&
 \TT_{\kappa,h}[\varphi]\le {\TT}_{\kappa,h}[\psi] -\frac{\theta_\kappa^2h^2}{2}e^{-\theta_\kappa h}\E[(\psi-\varphi)(t+h,\hat X_h^{t,x,\kappa})].
 \e*
In particular, if $\kappa_h$ satisfies \reff{thetakappah},
then 
 \b*
 \varphi\le\psi
 &\Longrightarrow&
 \TT_{\kappa_h,h}[\varphi]\le {\TT}_{\kappa_h,h}[\psi] +C h\E[(\psi-\varphi)(t+h,\hat X_h^{t,x,\kappa_h})]
 \e*
 for some constant $C$.
\end{Corollary}
\proof 
Let $\theta_\kappa$ be as in Lemma \ref{lemmonjumpkappa} and define $\varphi_\kappa(t,x):=e^{\theta_\kappa(T-t)}\varphi(t,x)$ and $\psi_\kappa(t,x):=e^{\theta_\kappa(T-t)}\psi(t,x)$. By Lemma \ref{lemmonjumpkappa}, 
\b*
\overline{\TT}_{\kappa,h}[\varphi_\kappa]&\le&\overline{\TT}_{\kappa,h}[\psi_\kappa].
\e*
By multiplying both sides by $e^{-\theta_\kappa(T-t)}$, we have
\begin{align*}
\left (e^{-\theta_\kappa h}(1+\theta_\kappa h)-1\right )&\E[\varphi(t+h,\hat X_h^{t,x,\kappa})]+\TT_{\kappa,h}[\varphi]\\
\le&\left (e^{-\theta_\kappa h}(1+\theta_\kappa h)-1\right )\E[\psi(t+h,\hat X_h^{t,x,\kappa})]+\TT_{\kappa,h}[\psi].
\end{align*}
So,
\b*
\TT_{\kappa,h}[\varphi]
&\le&
\left (e^{-\theta_\kappa h}(1+\theta_\kappa h)-1\right )\E[(\psi-\varphi)(t+h,\hat X_h^{t,x,\kappa})]+\TT_{\kappa,h}[\psi].
\e*
But, $e^{-\theta_\kappa h}(1+\theta_\kappa h)-1\le-\frac{\theta_\kappa^2h^2}{2}e^{-\theta_\kappa h}$. So,
\b*
\TT_{\kappa,h}[\varphi]
&\le&
-\frac{\theta_\kappa^2h^2}{2}e^{-\theta_\kappa h}\E[(\psi-\varphi)(t+h,\hat X_h^{t,x,\kappa})]+\TT_{\kappa,h}[\psi].
\e*
which \reff{thetakappah} provides the result. \ep
\vspace*{5mm}

In order to provide a uniform bound on $v^{\kappa,h}$, we bound $u^{\kappa,h}$ with respect to $\theta_\kappa$ as in the following Lemma.
\begin{Lemma}\label{lemstabjumpbar}
Let $\varphi$ and $\psi:[0,T]\times\R^d\longrightarrow\R$ be two $L^\infty-$bounded functions. Then 
 \b*
|\overline{\bf T}_{\kappa,h}[\varphi]-\overline{\bf T}_{\kappa,h}[\psi]|_\infty\le|\varphi-\psi|_\infty(1+(C+\theta_\kappa) h)
 \e*
where $C=\sup_{\alpha,\beta}|c^{\alpha,\beta}|_\infty$.
In particular, if $g$ is $L^\infty-$bounded, for a fixed $\kappa$ the family $(u^{\kappa,h}(t,\cdot))_h$ defined in \reff{schemejump} is $L^\infty-$bounded, uniformly in $h$ by 
\b*
(\overline{C}+|g|_\infty) e^{(C+\theta_\kappa)(T-t_i)}.
\e*
\end{Lemma}
\proof Let $f:=\varphi-\psi$. Then, by Assumption {\bf Inf-Sup} and the same argument as in the proof of Lemma \ref{lemmonjumpkappa}, 
\b*
\overline{\bf T}_{\kappa,h}[\varphi]-\overline{\bf T}_{\kappa,h}[\psi]
\!\!\!&\!\!\le\!\!&\!\!\!
\E\biggl[ f(\hat X_h)\biggl(1-a^{-1}\cdot a^{\alpha^*,\beta^*}+h \Bigl(|A^{\alpha^*,\beta^*}_h|^2 +c^{\alpha^*,\beta^*}+\theta_\kappa\\
\!\!\!&\!\!\!\!&\!\!\!-\int_{\{|z|\ge\kappa\}}\!\!\!\!\!\!\nu(dz)-\frac{1}{4}\Bigl(b^{\alpha^*,\beta^*}-\int_{\{1>|z|\ge\kappa\}}\!\!\!\!\!\!\!\!\!\!\!\!\!\!\!\!\eta^{\alpha^*,\beta^*}(z)\nu(dz)\Bigr )^{\rm T}({a^{\alpha^*,\beta^*}})^{-}\\
\!\!\!&\!\!\!\!&\!\!\!\times\Bigl(b^{\alpha^*,\beta^*}-\int_{\{1>|z|\ge\kappa\}}\!\!\!\!\!\!\!\!\!\!\!\!\!\!\!\!\eta^{\alpha^*,\beta^*}(z)\nu(dz)\Bigr)\Bigr)\biggr)+h\hat\nu_h^{\eta^{\alpha^*,\beta^*},1}(f)\biggr],
\e*
where $A^{\alpha^*,\beta^*}_h$ is given by \reff{Ah}.
On the other hand, 
\b*
\left|\hat\nu_h^{\eta^{\alpha^*,\beta^*},1}(f)\right|\le|f|_\infty\int_{\{|z|\ge\kappa\}}\!\!\!\!\!\!\!\!\!\!\!\!\nu(dz)
\e*
Therefore , 
\begin{align*}
 \overline{\bf T}_{\kappa,h}[\varphi]&-\overline{\bf T}_{\kappa,h}[\psi]
 \le 
 |f|_\infty\E\biggl[\Bigl|1-a^{-1}\cdot a^{\alpha^*,\beta^*}+h \Bigl(|A^{\alpha^*,\beta^*}_h|^2+c^{\alpha^*,\beta^*}+\theta_\kappa\\
\!\!\!&\!\!\!\!\!\!\!-\frac{1}{4}\Bigl(b^{\alpha^*,\beta^*}-\int_{\{1>|z|\ge\kappa\}}\!\!\!\!\!\!\!\!\!\!\!\!\!\!\!\!\eta^{\alpha^*,\beta^*}(z)\nu(dz)\Bigr)^{\rm T}({a^{\alpha^*,\beta^*}})^{-}
\Bigl(b^{\alpha^*,\beta^*}-\int_{\{1>|z|\ge\kappa\}}\!\!\!\!\!\!\!\!\!\!\!\!\!\!\!\!\eta^{\alpha^*,\beta^*}(z)\nu(dz)\Bigr)\Bigr)\Bigr|\biggr].
\end{align*}
By Assumption {\bf IHJB1} and \reff{thetakappa}, $1-a^{-1}\cdot a^{\alpha^*,\beta^*}$ and 
\begin{align*}
&c^{\alpha^*,\beta^*}+\theta_\kappa
-\frac{1}{4}\Bigl(b^{\alpha^*,\beta^*}-\int_{\{1>|z|\ge\kappa\}}\!\!\!\!\!\!\!\!\!\!\!\!\!\!\!\!\eta^{\alpha^*,\beta^*}(z)\nu(dz)\Bigr)^{\rm T}({a^{\alpha^*,\beta^*}})^{-}\Bigl(b^{\alpha^*,\beta^*}-\int_{\{1>|z|\ge\kappa\}}\!\!\!\!\!\!\!\!\!\!\!\!\!\!\!\!\eta^{\alpha^*,\beta^*}(z)\nu(dz)\Bigr)
\end{align*}
are positive. Therefore, one can write
\begin{align}\label{stabarg1}
\overline{\bf T}_{\kappa,h}[\varphi]&-\overline{\bf T}_{\kappa,h}[\psi]
 \le 
 |f|_\infty\biggl(1-a^{-1}\cdot a^{\alpha^*,\beta^*}+h\Bigl( \E[|A^{\alpha^*,\beta^*}_h|^2]+c^{\alpha^*,\beta^*}+\theta_\kappa\\
\!\!\!&\!\!\!\!\!\!\!-\frac{1}{4}\Bigl(b^{\alpha^*,\beta^*}-\int_{\{1>|z|\ge\kappa\}}\!\!\!\!\!\!\!\!\!\!\!\!\!\!\!\!\eta^{\alpha^*,\beta^*}(z)\nu(dz)\Bigr)^{\rm T}({a^{\alpha^*,\beta^*}})^{-}
\Bigl(b^{\alpha^*,\beta^*}-\int_{\{1>|z|\ge\kappa\}}\!\!\!\!\!\!\!\!\!\!\!\!\!\!\!\!\eta^{\alpha^*,\beta^*}(z)\nu(dz)\Bigr)\Bigr)\biggr). \nonumber
\end{align}
 But,  Notice that
\begin{align*}
 \E[|A^{\alpha^*,\beta^*}_h|^2]&=h^{-1}a^{-1}\cdot a^{\alpha^*,\beta^*}\\
&+\frac{1}{4}\Bigl(b^{\alpha^*,\beta^*}-\int_{\{1>|z|\ge\kappa\}}\!\!\!\!\!\!\!\!\!\!\!\!\!\!\!\!\eta^{\alpha^*,\beta^*}(z)\nu(dz)\Bigr )^{\rm T}{a^{\alpha^*,\beta^*}}^{-1}\Bigl(b^{\alpha^*,\beta^*}-\int_{\{1>|z|\ge\kappa\}}\!\!\!\!\!\!\!\!\!\!\!\!\!\!\!\!\eta^{\alpha^*,\beta^*}(z)\nu(dz)\Bigr).
\end{align*}
By replacing $\E[|A^{\alpha^*,\beta^*}_h|^2]$ into \reff{stabarg1}, one obtains
\b*
\overline{\bf T}_{\kappa,h}[\varphi]-\overline{\bf T}_{\kappa,h}[\psi]
 &\le& 
 |f|_\infty(1+h(c^{\alpha^*,\beta^*}+\theta_\kappa))\\
&\le& |f|_\infty(1+(C+\theta_\kappa)h),
 \e*
with $C=\sup_{\alpha,\beta}|c^{\alpha,\beta}|_\infty$. By changing the role of $\varphi$ and $\psi$ and implementing the same argument, one obtains
\b*
 \left|\overline{\bf T}_{\kappa,h}[\varphi]-\overline{\bf T}_{\kappa,h}[\psi]\right|_\infty
&\le& |f|_\infty(1+(C+\theta_\kappa)h).
 \e*
\\

To prove that the family $(u^{\kappa,h})_h$ is bounded, we proceed by backward induction as in Lemma 3.14 in \cite{ftw}. By choosing in the first part of the proof $\varphi\equiv \bar u^{\kappa,h}(t_{i+1},.)$ and $\psi\equiv 0$, we see that
\b*
|u^{\kappa,h}(t_i,\cdot)|_\infty
&\le&
h\overline{C}e^{\theta_\kappa(T-t_i)}
+
|u^{\kappa,h}(t_{i+1},\cdot)|_\infty(1+(C+\theta_\kappa) h),
 \e*
where $\overline{C}:=\sup_{\alpha,\beta}|k^{\alpha,\beta}|_\infty$. It follows from the discrete Gronwall inequality that 
\b*
|u^{\kappa,h}(t_i,\cdot)|_\infty
\le (\overline{C}(T-t_i)+|g|_\infty) e^{(C+\theta_\kappa)(T-t_i)}.
\e*
\ep\\

Define 
\be\label{vbarukappah}
\bar v^{\kappa,h}:=e^{-\theta_\kappa(T-t)}u^{\kappa,h}.
\ee
 Next Corollary provides a bound for $v^{\kappa,h}$ uniformly on $\kappa$ and $h$.
\begin{Corollary}\label{corbarvkappah}
$\bar v^{\kappa,h}$ is bounded uniformly on $h$ and $\kappa$, and
\b*
|v^{\kappa,h}-\bar v^{\kappa,h}|_\infty&\le&K\theta_\kappa^2 h~~~\text{for some constant}~~~K.
\e*
If also, $\kappa_h$ satisfies \reff{thetakappah}, 
then
\b*
\lim_{h\to0}|v^{\kappa_h,h}-\bar v^{\kappa_h,h}|_\infty&=&0.
\e*
\end{Corollary}
\proof By Lemma \ref{lemstabjumpbar} for fixed $\kappa$, we have:
\b*
|u^{\kappa,h}(t,.)|_\infty
&\le&
(\overline{C}+|g|_\infty) e^{(C+\theta_\kappa)(T-t)}.
\e*
Therefore, 
\b*
|\bar v^{\kappa,h}(t,.)|_\infty
&\le&
(\overline{C}+|g|_\infty) e^{C(T-t)}.
\e*

For the next part, define $\bar u^{\kappa,h}(t,x)=e^{\theta_\kappa(T-t)}v^{\kappa,h}(t,x)$. Direct calculations shows that
\b*
\bar u^{\kappa,h}=e^{\theta_\kappa h}(1-\theta_\kappa h)\E\left[\bar u^{\kappa,h}\left (t+h,\hat X_h^{t,x,\kappa}\right)\right]
                   +h\overline{F}_{\kappa,h}\left(t,x,\Dc_h\bar u^{\kappa,h},\bar u^{\kappa,h}(t+h,\cdot)\right).
\e*
By an argument similar to Lemma 3.19 in \cite{ftw}, we have
\be\label{ubaru}
|(u^{\kappa,h}-\bar u^{\kappa,h})(t,\cdot)|_\infty 
& \le &
\frac12\theta_\kappa^2h^2|\bar u^{\kappa,h}(t+h,\cdot)|_\infty\\
&&+\nonumber(1+(C+\theta_\kappa)h)|(u^{\kappa,h}-\bar u^{\kappa,h})(t+h,\cdot)|_\infty,
\ee
where $C$ is as in Lemma \ref{lemstabjumpbar}.
By repeating the proof of Lemma \ref{lemstabjumpbar} for $\bar u^{\kappa,h}$, one can conclude,
\b*
|\bar u^{\kappa,h}(t,\cdot)|_\infty
&\le&
(\overline{C}+|g|_\infty) e^{(C+\theta_\kappa)(T-t)}(1+\frac{\theta_\kappa h}{2}).
\e*
So, by multiplying \ref{ubaru} by $e^{\theta_\kappa(T-t)}$, we have
\b*
|(\bar v^{\kappa,h}-v^{\kappa,h})(t,\cdot)|_\infty&\le&\frac12\tilde C\theta_\kappa^2h^2e^{C(T-t)}(1+\frac{\theta_\kappa h}{2})e^{-\theta_\kappa h}\\
&&+e^{-\theta_\kappa h}(1+(C+\theta_\kappa)h)|(\bar v^{\kappa,h}-v^{\kappa,h})(t+h,\cdot)|_\infty,
\e*
for some constant $\tilde C$.
Because $e^{-\theta_\kappa h}(1+(C+\theta_\kappa)h)\le e^{Ch}$, one can deduce from discrete Gronwall inequality that
\b*
|(\bar v^{\kappa,h}-v^{\kappa,h})(t,\cdot)|_\infty 
&\le&
K\theta_\kappa^2h,
\e*
for some constant $K$ independent of $\kappa$ which provides the second part of the theorem.\ep
\vspace*{5mm}

We continue with the following consistency Lemma.
\begin{Lemma}\label{lemconsjump}
Let $\varphi$ be a smooth function with the bounded derivatives. Then for all $(t,x)\in[0,T]\times\R^d$:
 \b*
 \lim_{\begin{array}{c}
\vspace*{-2mm}
             {\scriptscriptstyle (t',x')\to(t,x)}
             \\
\vspace*{-2mm}
             {\scriptscriptstyle (h,c)\to (0,0)}
             \\
             {\scriptscriptstyle t'+h\le T}
             \end{array}}
\!\!\!\!\!\!\!\!\!\!\!\! \frac{\varphi(t',x')-\TT_{\kappa,h}[c+\varphi](t',x')}{h}
 \!\!\!&=&\!\!\!
 -\left(\mathcal{L}^X\varphi+F(\cdot,\varphi,D\varphi,D^2\varphi,\varphi(t,\cdot))\right)(t,x).
 \e*
\end{Lemma}
\proof The proof is straightforward by Lebesgue dominated convergence Theorem.\ep\\
\vspace*{5mm}

To complete the convergence argument, we need to proof the the approximate solution $v^{\kappa_h,h}$ converge to  the final condition as 
\begin{Lemma}\label{lemregvkappah}
Let  $\kappa_h$ satisfy \reff{thetakappah}, then
$\bar v^{\kappa_h,h}$ is uniformly Lipschitz with respect to $x$.
\end{Lemma}
\proof We report the following calculation in the one-dimensional case $d=1$ in order to simplify the presentation. \\
For fixed $t\in[0,T-h]$, we argue as in the proof of Lemma \ref{lemstabjumpbar} to see that for $x,x'\in\R$ with $x> x'$:
\b*
 u^{\kappa,h}(t,x)-u^{\kappa,h}(t,x')
 \!\!\!&=&\!\!\!
\E\Bigl[\left(u^{\kappa,h}(t+h,\hat X^{t,x})-u^{\kappa,h}(t+h,\hat X^{t,x'})\right)]\\
 &&+h\left(\inf_{\alpha}\sup_{\beta}\bar\Jc_\kappa^ {\alpha,\beta}[\what{u^{\kappa,h}}](t+h,x)-
 \inf_{\alpha}\sup_{\beta}\bar\Jc_\kappa^ {\alpha,\beta}[\what{u^{\kappa,h}}](t+h,x')\right)\\
&\le&
\E\Bigl[\left(u^{\kappa,h}(t+h,\hat X^{t,x})-u^{\kappa,h}(t+h,\hat X^{t,x'})\right)]\\
 &&+h\left(\sup_{\beta}\bar\Jc_\kappa^ {\alpha,\beta}[\what{u^{\kappa,h}}](t+h,x)-
 \inf_{\alpha}\bar\Jc_\kappa^ {\alpha,\beta}[\what{u^{\kappa,h}}](t+h,x')\right).
\e*
Observe that by \reff{Phi}, one can write
\b*
u^{\kappa,h}(t,x)-u^{\kappa,h}(t,x')
 &\le&
\E\Bigl[\left(u^{\kappa,h}(t+h,\hat X^{t,x})-u^{\kappa,h}(t+h,\hat X^{t,x'})\right)\Bigr]\\
 &&+h\Bigl(\bar\Phi^{\alpha,\beta}[\what{u^{\kappa,h}},\what{u^{\kappa,h}}](t+h,x,x')\Bigr) ,
\e*
where $\bar \Phi$ is defined in the proof of Lemma \ref{lemstabjumpbar}.
By Assumption {\bf Inf-Sup}, there exists $(\alpha^*,\beta^*)$  such that
\b*
\bar\Phi^{\alpha^*,\beta^*}[\what{u^{\kappa,h}},\what{u^{\kappa,h}}](t+h,x,x')=\bar\Jc_\kappa^ {\alpha^*,\beta^*}[\what{u^{\kappa,h}}](t+h,x)-\bar\Jc_\kappa^ {\alpha^*,\beta^*}[\what{u^{\kappa,h}}](t+h,x').
\e*
Therefore,
\b*
u^{\kappa,h}(t,x)-u^{\kappa,h}(t,x')
 \!\!\!&\le&\!\!\!
\E\Bigl[\left(u^{\kappa,h}(t+h,\hat X^{t,x})-u^{\kappa,h}(t+h,\hat X^{t,x'})\right)]\\
  &&+h\Bigl(
\bar\Jc_\kappa^*[\what{u^{\kappa,h}}](t+h,x)- \bar\Jc_\kappa^*[\what{u^{\kappa,h}}](t+h,x')
\Bigr).
\e*

For the other in equality we do the same except that when we 
\b*
 u^{\kappa,h}(t,x)-u^{\kappa,h}(t,x')
 &\le&
A+hB+hC,
\e*
where 
\b*
A
&:=&
\E\Bigl[\left(u^{\kappa,h}(t+h,\hat X^{t,x})-u^{\kappa,h}(t+h,\hat X^{t,x'})\right)]\\
&&+h\left(\bar\Jc_\kappa^ {\alpha^*,\beta^*}[\what{u^{\kappa,h}}](t+h,x)-
 \bar\Jc_\kappa^ {\alpha^*,\beta^*}[\what{\tilde u^{\kappa,h}}](t+h,x)\right),
\e*
with $\tilde u^{\kappa,h}(y)=u^{\kappa,h}(y+x'-x)$,
\b*
B
&:=&
\bar\Jc_\kappa^ {\alpha^*,\beta^*}[\what{\tilde u^{\kappa,h}}](t+h,x)-
\bar\Jc_\kappa^ {\alpha^*,\beta^*}[\what{u^{\kappa,h}}](t+h,x'),
\e*
and
\b*
C:=\hat\nu_h^{\alpha^*,\beta^*,1}(u^{\kappa,h}(t+h,\cdot))(x)-\hat\nu_h^{\alpha^*,\beta^*,1}(u^{\kappa,h}(t+h,\cdot))(x').
\e*
We continue the proof in the following steps. \\

\no{\bf Step 1.} 
\b*
C=h^{-1}\E\Bigl[\Bigl(u^{\kappa,h}(t+h,\Xh^{*,x})-u^{\kappa,h}(t+h,\Xh^{*,x'})\Bigr)N_h^{\kappa}\Bigr],
\e*
where $\Xh^{*,x}:=x+\sum_{i=1}^{N_h^{\kappa}}\eta^{\alpha^*,\beta^*}(x,Z_i)$ with $Z_i$s are i.i.d. random variables distributed as $\frac{\nu(dz)}{\lambda_\kappa}$.

\no{\bf Step 2.} 
By the definition of $\bar\Jc_\kappa^ {\alpha,\beta}$,
\b*
B
&=&
\frac12(a^ {\alpha^*,\beta^*}(x)-a^ {\alpha^*,\beta^*}(x'))\Dc_h^2u^{\kappa,h}(t+h,x')
+(b_\kappa^ {\alpha^*,\beta^*}(x)-b_\kappa^ {\alpha^*,\beta^*}(x'))\Dc_h^1u^{\kappa,h}(t+h,x')\\
&&\hspace*{2cm}+(c^ {\alpha^*,\beta^*}(x)-c^ {\alpha^*,\beta^*}(x'))\Dc_h^0u^{\kappa,h}(t+h,x')+k^ {\alpha^*,\beta^*}(x)-k^ {\alpha^*,\beta^*}(x'),
\e*
where $b^{\alpha,\beta}_\kappa(x):=b^{\alpha,\beta}(x)-\int_{\{1>|z|\ge\kappa\}}\eta^{\alpha,\beta}(x,z)\nu(dz)$. 
On the other hand,
\b*
\Dc_h^k=\E\left[Du^{\kappa,h}(t+h,\Xh_h^{x'})\left(\frac{W_h}{h}\sigma^{-1}(x')\right)^{k-1}\right],~~\text{for}~~k=1,2.
\e*
So,
\b*
B
&\le&
\E\biggl[
\frac12(a^ {\alpha^*,\beta^*}(x)-a^ {\alpha^*,\beta^*}(x'))Du^{\kappa,h}(t+h,\Xh_h^{x'})\frac{W_h}{h}\sigma^{-1}(x')\\
&&+(b_\kappa^ {\alpha^*,\beta^*}(x)-b_\kappa^ {\alpha^*,\beta^*}(x'))Du^{\kappa,h}(t+h,\Xh_h^{x'})+(c^ {\alpha^*,\beta^*}(x)-c^ {\alpha^*,\beta^*}(x'))u^{\kappa,h}(t+h,\Xh_h^{x'})
\biggr]\\
&&+f^ {\alpha^*,\beta^*}(x)-f^ {\alpha^*,\beta^*}(x').
\e*

\no{\bf Step 3.} 
By the definition of $\bar\Jc_\kappa^ {\alpha,\beta}$, one can observe that
\begin{align*}
\bar\Jc_\kappa^ {\alpha^*,\beta^*}[u^{\kappa,h}](t+h,x)&-
 \bar\Jc_\kappa^ {\alpha^*,\beta^*}[\tilde u^{\kappa,h}](t+h,x)\\
&=\frac12 a^{\alpha^*,\beta^*}(x)\delta^{(2)}+b_\kappa^*(x)\delta{(1)} +c_\kappa^*(x)\delta{(0)}
\end{align*}
where $c_\kappa^*$ and $b_\kappa^*$ are defind in the proof of Lemma \ref{lemmonjumpkappa}, and
\b*
\delta^{(k)}
&=&
\E\left[D^ku^{\kappa,h}(t+h,\Xh_h^x)-D^ku^{\kappa,h}(t+h,\Xh_h^{x'})\right]~~\text{for}~~k=0,1,2.
\e*
By Lemma 2.1 in \cite{ftw},  for $k=1$ and $2$
\b*
\delta^{(k)}
&=&
\E\Bigl[\left(u^{\kappa,h}(t+h,\Xh_h^x)-u^{\kappa,h}(t+h,\Xh_h^{x'})\right)H_h^k(t,x)\\
&&+u^{\kappa,h}(t+h,\Xh_h^{x'})H_h^k(t,x)\left(1-\frac{\sigma^{k}(x)}{\sigma^{k}(x')}\right)\Bigr]\\
&=&
\E\Bigl[\left(u^{\kappa,h}(t+h,\Xh_h^x)-u^{\kappa,h}(t+h,\Xh_h^{x'})\right)H_h^k(t,x)\\
&&+Du^{\kappa,h}(t+h,\Xh_h^{x'})\left(\frac{W_h}{h}\right)^{k-1}\sigma(x')\left(\sigma^{-k}(x)-\sigma^{-k}(x')\right)\Bigr].
\e*
Therefore, one can write
\b*
A
&\le&
\E\Bigl[\left(u^{\kappa,h}(t+h,\Xh_h^x)-u^{\kappa,h}(t+h,\Xh_h^{x'})\right)\\
&&\hspace*{2cm}\times\Bigl(1-\bar a^*+\bar a^*N^2
+hc_\kappa^*+b_\kappa^*N\sqrt{h}\Bigr)(x)\\
&&\hspace*{2cm}+hb_\kappa^*(x')Du^{\kappa,h}(t+h,\Xh_h^{x'})\sigma(x')\left(\sigma^{-1}(x)-\sigma^{-1}(x')\right)\\
&&\hspace*{2cm}+a^*(x')Du^{\kappa,h}(t+h,\Xh_h^{x'})\sqrt{h}N\sigma(x')\left(\sigma^{-2}(x)-\sigma^{-2}(x')\right)\Bigr],
\e*
where $a^*:=\frac12a^{\alpha^*,\beta^*}$, $\bar a^*:=\frac12a^{-1}a^{\alpha^*,\beta^*}$, $c^*:=c^{\alpha^*,\beta^*}$, $c_\kappa^*:=c^*+\theta_\kappa$, and $b_\kappa^*:=b_\kappa^{\alpha^*,\beta^*}$.

\no{\bf Step 4.} By dividing both sides by $x-x'$ and taking the limit we have: 
\b*
Du^{\kappa,h}(t,x)
&\le&
\E\Bigl[
Du^{\kappa,h}(t+h,\Xh_h^{x})\Biggl(\left(1+h\tilde\mu'_\kappa+\sqrt{h}\sigma'N+\Jt_{\kappa,h}\right)\\
&&\times\Bigl(1-\bar a^*+\bar a^*N^2
+hc_\kappa^*+b_\kappa^*N\sqrt{h}\Bigr)\\
&&+h\Bigl((b_\kappa^*)'-b_\kappa^*\frac{\sigma'}{\sigma}\Bigr)+\Bigl(\frac12(a^{\alpha^*,\beta^*})'\sigma^{-1}-a^{\alpha^*,\beta^*}\frac{\sigma'}{\sigma^2}\Bigr)\sqrt{h}N
\biggr)\\
&&+Du^{\kappa,h}(t+h,\Xh_h^{*,x})\left(1+\mu^*h+\Jt_{\kappa,h}^{'*}\right)N_h^\kappa\Bigr]+Ce^{\theta\kappa(T-t)}h,
\e*
where $\Jt_{\kappa,h}:=\int_{\{|z|>\kappa\}}\eta(z)\Jt([0,h],dz)$, $\Jt^{'*}_{\kappa,h}:=\int_{\{|z|>\kappa\}}\eta'(z)\Jt([0,h],dz)$, and $N_h^\kappa$ is a Poisson process with intensity $\lambda_\kappa:=\int_{\{|z|>\kappa\}}\nu(dz)$. 

\no Let $L_{t}:=|Du^{\kappa,h}(t,\cdot)|_\infty$. Then
\b*
\E\Bigl[Du^{\kappa,h}(t+h,\Xh_h^{*,x})\left(1+\mu^*h+\Jt_{\kappa,h}^{'*}\right)N_h^\kappa\Bigr]
&\le&
L_{t+h}Ch\left(\lambda_\kappa+\lambda'^{*}_\kappa\right),
\e*
where $\lambda'^{*}_\kappa:=\int_{\{|z|>\kappa\}}\eta'^{*}(z)\nu(dz)$.
Let $G:=N+\frac{b_\kappa^*\sigma}{2}\sqrt{h}$. By the change of measure 
\b*
\frac{d\Q}{d\P}:=\exp\left(-\frac{(b_\kappa^*\sigma)^2}{4}h+\frac{b_\kappa^*\sigma}{2}\sqrt{h}N\right),
\e*
we have $G\sim \Nc(0,1)$ under $\Q$ and one can write
\b*
Du^{\kappa,h}(t,x)
&\le&
\E^{\Q}\Bigl[\frac{d\P}{d\Q}
Du^{\kappa,h}(t+h,\Xh_h^{x})\Bigl(\Bigl(1+h(\tilde\mu'_\kappa-\frac{b_\kappa^*\sigma}{2})+\sqrt{h}\sigma'G+\Jt_{\kappa,h}\Bigr)\\
&&\times\Bigl(1-\bar a^*+\bar a^*G^2
+h(c_\kappa^*-\frac{(b_\kappa^*\sigma)^2}{2})\Bigr)\\
&&+h\Bigl((b_\kappa^*)'-b_\kappa^*\frac{\sigma'}{\sigma}-\frac{b_\kappa^*\sigma}{2}\Bigr)+\Bigl(\frac12(a^{\alpha^*,\beta^*})'\sigma^{-1}-a^{\alpha^*,\beta^*}\frac{\sigma'}{\sigma^2}\Bigr)\sqrt{h}G
\Bigr)
\Bigr]\\
&&+L_{t+h}Ch\left(\lambda_\kappa+\lambda'^{*}_\kappa\right)+Ce^{\theta\kappa(T-t)}h,
\e*

%
\no{\bf Step 4.} Notice that $1-\bar a^*+a^*G^2+h(c_\kappa^*-\frac{(b_\kappa^*\sigma)^2}{2})$ is positive and therefore, one can take $\frac{Z}{\E^\Q[Z]}$ as a density for the new measure $\Q^Z$. So, 
\b*
Du^{\kappa,h}(t,x)
&\le&
\E^{\Q^Z}\Bigl[\frac{d\P}{d\Q}
Du^{\kappa,h}(t+h,\Xh_h^{x})\Bigl(\Bigl(1+h(\tilde\mu'_\kappa-\frac{b_\kappa^*\sigma}{2})+\sqrt{h}\sigma'G+\Jt_{\kappa,h}\Bigr)\\
&&+Z^{-1}\Bigl(h\Bigl((b_\kappa^*)'-b_\kappa^*\frac{\sigma'}{\sigma}-\frac{b_\kappa^*\sigma}{2}\Bigr)+\Bigl(\frac12(a^{\alpha^*,\beta^*})'\sigma^{-1}-a^{\alpha^*,\beta^*}\frac{\sigma'}{\sigma^2}\Bigr)\sqrt{h}G\Bigr)
\Bigr)
\Bigr]\\
&&+L_{t+h}Ch\left(\lambda_\kappa+\lambda'^{*}_\kappa\right)+Ce^{\theta\kappa(T-t)}h.
\e*
So,
\b*
Du^{\kappa,h}(t,x)
&\le&
\E^{\Q^Z}\Bigl[\left(\frac{d\P}{d\Q}\right)^2
(Du^{\kappa,h}(t+h,\Xh_h^{x}))^2\Bigr]^\frac12\E^{\Q^Z}\Bigl[\Bigl(\Bigl(1+h(\tilde\mu'_\kappa-\frac{b_\kappa^*\sigma}{2})+\sqrt{h}\sigma'G+\Jt_{\kappa,h}\Bigr)\\
&&+Z^{-1}\Bigl(h\Bigl((b_\kappa^*)'-b_\kappa^*\frac{\sigma'}{\sigma}-\frac{b_\kappa^*\sigma}{2}\Bigr)+\Bigl(\frac12(a^{\alpha^*,\beta^*})'\sigma^{-1}-a^{\alpha^*,\beta^*}\frac{\sigma'}{\sigma^2}\Bigr)\sqrt{h}G\Bigr)
\Bigr)^2
\Bigr]^\frac12\\
&&+L_{t+h}Ch\left(\lambda_\kappa+\lambda'^{*}_\kappa\right)+Ce^{\theta\kappa(T-t)}h.
\e*
Notice that
\b*
\E^{\Q^Z}\Bigl[\left(\frac{d\Q}{d\P}\right)^2
(Du^{\kappa,h}(t+h,\Xh_h^{x}))^2\Bigr]
&\le&
L_{t+h}^2\exp(\frac{1}{4}(b_\kappa^*\sigma)^2h).
\e*
On the other hand,
\begin{align*}
&\E^{\Q^Z}\Bigl[\frac{d\Q}{d\P}\Bigl(\Bigl(1+h(\tilde\mu'_\kappa-\frac{b_\kappa^*\sigma}{2})+\sqrt{h}\sigma'G+\Jt_{\kappa,h}\Bigr)+Z^{-1}\Bigl(h\Bigl((b_\kappa^*)'-b_\kappa^*\frac{\sigma'}{\sigma}-\frac{b_\kappa^*\sigma}{2}\Bigr)\\
&\hspace*{3cm}+\Bigl(\frac12(a^{\alpha^*,\beta^*})'\sigma^{-1}-a^{\alpha^*,\beta^*}\frac{\sigma'}{\sigma^2}\Bigr)\sqrt{h}G\Bigr)
\Bigr)^2\Bigr]\\
&=
\E\Bigl[Z\Bigl(\Bigl(1+h(\tilde\mu'_\kappa-\frac{b_\kappa^*\sigma}{2})+\sqrt{h}\sigma'G+\Jt_{\kappa,h}\Bigr)+Z^{-1}\Bigl(h\Bigl((b_\kappa^*)'-b_\kappa^*\frac{\sigma'}{\sigma}-\frac{b_\kappa^*\sigma}{2}\Bigr)\\
&\hspace*{3cm}+\Bigl(\frac12(a^{\alpha^*,\beta^*})'\sigma^{-1}-a^{\alpha^*,\beta^*}\frac{\sigma'}{\sigma^2}\Bigr)\sqrt{h}G\Bigr)
\Bigr)^2\Bigr].
\end{align*}
By calculation of the right hand side of the above equality, one can observe that all the terms of order $\sqrt{h}$ vanish and we have:
\begin{align*}
&\E^{\Q^Z}\Bigl[\frac{d\Q}{d\P}\Bigl(\Bigl(1+h(\tilde\mu'_\kappa-\frac{b_\kappa^*\sigma}{2})+\sqrt{h}\sigma'G+\Jt_{\kappa,h}\Bigr)+Z^{-1}\Bigl(h\Bigl((b_\kappa^*)'-b_\kappa^*\frac{\sigma'}{\sigma}-\frac{b_\kappa^*\sigma}{2}\Bigr)\\
&\hspace*{3cm}+\Bigl(\frac12(a^{\alpha^*,\beta^*})'\sigma^{-1}-a^{\alpha^*,\beta^*}\frac{\sigma'}{\sigma^2}\Bigr)\sqrt{h}G\Bigr)
\Bigr)^2\Bigr]^\frac12\\
&\le
\Bigl(1+h\Bigl(c^*+\theta_\kappa-\frac{(b^*_\kappa)^2}{4a^*}-b_\kappa^*\sigma\sigma'+(b^*_\kappa)'-\frac{b_\kappa^*\sigma'}{\sigma}-\frac{b_\kappa^*\sigma}{2}+O(h\theta_\kappa^2)\Bigr)\Bigr)^\frac12.
\end{align*}
Therefore,
by the choice of $\kappa_h$, for $h$ small enough we have
\b*
L_t
&\le& L_{t+h}\exp\Bigl(\frac12h(C+\theta{\kappa_h}-b{\kappa_h}^*\sigma\sigma'+(b^*{\kappa_h})'-\frac{b{\kappa_h}^*\sigma'}{\sigma}-\frac{b_{\kappa_h}^*\sigma}{2}+2\lambda_{\kappa_h}+2\lambda'^*{\kappa_h})\Bigl)
+Ce^{\theta_{\kappa_h}(T-t)}h\\
&\le&
L_{t+h}\exp\Bigl(h(C+\theta{\kappa_h})\Bigl)
+Ce^{\theta_{\kappa_h}(T-t)}h\\.
\e*

By discrete Gronwall inequality,
\b*
L_t
&\le&
(|Dg|_\infty+C(T-t))e^{(\theta_{\kappa_h}+C)(T-t)}.
\e*
Therefore by definition of $\bar v^{\kappa,h}$, we have
\b*
|D\bar v^{\kappa_h,h}|_1
&\le&
e^{C(T-t)}(|Dg|_\infty+C(T-t)).
\e* 
\ep

\begin{Lemma}\label{lemfinalcond}
Let $\kappa_h$ satisfies \reff{thetakappah}, then
\b*
\lim_{t\to T}\bar v^{\kappa,h}(t,x)
&=&
g(x).
\e*
\end{Lemma}
\proof
We follow the same notations as in the proof of the previous Lemma and  write
\b*
 u^{\kappa,h}(t,x)
 \!\!\!&=&\!\!\!
\E\Bigl[u^{\kappa,h}(t+h,\hat X^{t,x})\Bigr]+h\inf_{\alpha}\sup_{\beta}\bar\Jc_\kappa^ {\alpha,\beta}[\what{u^{\kappa,h}}](t+h,x)\\
&\le&
\E\Bigl[u^{\kappa,h}(t+h,\hat X^{t,x})\Bigr]+h\sup_{\beta}\bar\Jc_\kappa^ {\alpha,\beta}[\what{u^{\kappa,h}}](t+h,x).
\e*
Observe that by \reff{Phi}, one can write
\b*
u^{\kappa,h}(t,x)
 &\le&
\E\Bigl[u^{\kappa,h}(t+h,\hat X^{t,x})\Bigr]+h\Bigl(\bar\Phi^{\alpha,\beta}[\what{u^{\kappa,h}},0](t+h,x,x')\Bigr)+h\sup_{\alpha,\beta}|f^{\alpha,\beta}|_\infty,
\e*
By Assumption {\bf Inf-Sup}, there exists $(\alpha^*,\beta^*)$ so that
\b*
 u^{\kappa,h}(t,x)
 \!\!\!&\le&\!\!\!
\E\Bigl[u^{\kappa,h}(t+h,\hat X^{t,x})\Bigr]+h\bar\Jc_\kappa^ {\alpha^*,\beta^*}[\what{u^{\kappa,h}}](t+h,x)+h\bar C,
\e*
where $\bar C:=\sup_{\alpha,\beta}|f^{\alpha,\beta}|_\infty$. Therefore,  for any $j=i,\cdots,n-1$ one can write 
\b*
 u^{\kappa,h}(t_j,\Xh_{t_j}^{t_i,x})
 \!\!\!&\le&\!\!\!
\E^\Q_{t_j}\Bigl[u^{\kappa,h}(t_{j+1},\Xh_{t_{j+1}}^{t_i,x})\Bigl(1-a^*_j+a^*_jG_j^2
+hC^*_j\Bigr)\Bigr]+h\bar C.
\e*
where $a^*_j:=\bar a^*(t_j,\Xh_{t_j}^{t_i,x})$, $C^*_j:=(c_\kappa^*-\frac{(b_\kappa^*\sigma)^2}{2})(t_j,\Xh_{t_j}^{t_i,x})$ and $G_j$s are independent standard Gaussian random variables under the new equivalent measure $\Q$.
By the consecutive use of the above inequality and the fact that $1-a^*_j+a^*_jG^2
+hC^*_j$ is positive, one can write
\b*
 u^{\kappa,h}(t_i,x)
 \!\!\!&\le&\!\!\!
\E^\Q\Bigl[g(\Xh_T^{t_i,x})\prod_{j=i}^{n-1}\Bigl(1-a^*_j+a^*_jG^2
+hC_j^*\Bigr)\Bigr]+\bar C h\sum_{j=i}^{n-1}e^{\theta_\kappa t_j}.
\e*
Notice that in the above inequality  we used the fact that 
\b*
\E^\Q_{t_j}\Bigl[1-a^*_j+a^*_jG_j^2
+hC_j^*\Bigr]
&=&
1+h\E^\Q_{t_j}[C_j^*]\le 1+\theta_\kappa h.
\e*
On the other hand, $Z:=\prod_{j=i}^{n-1}\Bigl(1-a^*_j+a^*_jG_j^2
+hC^*_j\Bigr)$ is positive there for $frac{Z}{E^\Q[Z]}$ could be considered as a density of a new measure $\Q^Z$ with respect to $\P$. Therefore,
\b*
u^{\kappa,h}(t_i,x)
 \!\!\!&\le&\!\!\!
\E^\Q[Z]\E^{\Q^Z}\Bigl[g(\Xh_T^{t_i,x})\Bigr]+\bar C h\sum_{j=i}^{n-1}e^{\theta_\kappa t_j}.
\e*
By the definition of $\bar v^{\kappa,h}$, one can write
\b*
\bar v^{\kappa,h}(t_i,x)
 \!\!\!&\le&\!\!\!
e^{-\theta_\kappa (T-t_i)}\E^\Q[Z]\E^{\Q^Z}\Bigl[g(\Xh_T^{t_i,x})\Bigr]+e^{-\theta_\kappa (T-t_i)}\bar C h\sum_{j=i}^{n-1}e^{\theta_\kappa t_j}.
\e*
Therefore,
\b*
\bar v^{\kappa,h}(t_i,x)-g(x)
 \!\!\!&\le&\!\!\!
e^{-\theta_\kappa (T-t_i)}\E^\Q[Z]\E^{\Q^Z}\Bigl[|g(\Xh_T^{t_i,x})-g(x)|\Bigr]+C|g(x)|(T-t_i)+e^{-\theta_\kappa (T-t_i)}\bar C (T-t_i).
\e*
Notice that $g(\Xh_T^{t_i,x})-g(x)$ converges to zero $\P$-a.s. and therefore $\Q^Z$ a.s. as $(t_i,h)\to (T,0)$. So, by Lebesgue dominated convergence Theorem, 
\b*
\limsup_{(t_i,h)\to (T,0)} \bar v^{\kappa,h}(t_i,x)-g(x) \le 0.
\e*
By the similar argument one  can prove that:
\b*
\liminf_{(t_i,h)\to (T,0)} \bar v^{\kappa,h}(t_i,x)-g(x) \ge 0,
\e*
which compelets the proof.\ep\\

\begin{Remark}\label{remtholder}
{\rm
By extending the above proof as in the Lemma 3.17 and Corollary 3.18 of \cite{ftw}, one can proof that
\b*
|\bar v^{\kappa,h}(t,x)-g(x)|
&\le&
C(T-t)^\frac12.
\e*
Also, observe that by the similar argumnet as in \cite{ftw}, $v^{\kappa_h,h}$ is $\frac12$-H\"older on $t$ uniformly on $h$ and $x$.
}
\end{Remark}

So, the approximate solution $\bar v^{\kappa_h,h}$ both satisfies the requirement of the convergence established in \cite{barlessouganidis} and converges to a function $v$ locally uniformly. Moreover, $v$ is the unique viscosity solution of \reff{equationnl}-\reff{terminalnl}. So, by Corollary \ref{corbarvkappah}, the same assertion is true for $v^{\kappa,h}$.

\subsection{Rate of Convergence}
\label{subsecrate}
For local PDE's, the rate of convergence of probabilistic numerical scheme relies on the approximation of the solution of PDE by regular sub and super-solutions, the consistency for the scheme and the comparison principle derived from strict monotonicity. One can approximate the solution of the local PDE by a regular sub-solution and a almost regular super-solution from up and down, respectively. These approximations are provided by a switching system and Krylov method of shaking coefficients. Next, we  use the consistency Lemma 3.22 in \cite{ftw} to produce inequalities for the regular approximations plugged into the scheme. Then by comparison principle; Proposition 3.20 in \cite{ftw}; we obtain the bounds for the difference of approximate solution derived from scheme and regular approximate solution obtained from Krylov method and switching system. 

We continue this Subsection with establishing the same line of argument as in \cite{ftw} for nonlocal case. The generalization of the method we used in \cite{ftw} for the rate of convergence, is developed in \cite{biswasjakobsenkarlsen1} where the scheme needs to be consistent and satisfies comparison principle. Before, providing consistency and comparison principle result for the scheme \reff{overlinescheme}, we show that truncation error could be handled by the Theorem of continuous dependence for \reff{equationnl}-\reff{terminalnl}.
More precisely, if $v$ and $v^\kappa$ are solutions of \reff{equationnl}-\reff{terminalnl} and \reff{equationnlkappa}-\reff{terminalnlkappa}, respectively; then by Theorem 5.1 in \cite{biswasjakobsenkarlsen2}
\b*
|v-v^\kappa|_\infty&\le&C\left (\int_{0<|z|<\kappa}|z|^2\nu(dz)\right )^\frac12.
\e*
Therefore, By choosing $\kappa_h$ so that $\int_{0<|z|<\kappa_h}|z|^2\nu(dz)\le Ch^\frac12$, one can just concentrate on the rate of convergence of $v^{\kappa,h}$ to $v^{\kappa}$.

We shift to $\bar v^{\kappa_h,h}$ which is is derived from the strictly monotone scheme \reff{overlinescheme} and find the rate of convergence for $\bar v^{\kappa_h,h}$. The following Corollary shows that this shift do not effect the rate of convergence.

\begin{Corollary}\label{corbarvkappahrate}
For $F$ which satisfies {\bf IHJB}, $F(t,x,0,0,0,0)=0$. 
Then,
\b*
|\bar v^{\kappa_h,h}-v^{\kappa_h,h}|\le Ch\theta_{\kappa_h}^2.
\e*
In addition, if $\kappa_h$ is such that
\b*
\limsup_{h\to0}h^\frac{3}{4}\theta_{\kappa_h}^2<\infty,
\e*
then 
\b*
|\bar v^{\kappa_h,h}-v^{\kappa_h,h}|\le Ch^\frac{1}{4}
\e*
\end{Corollary}
\proof The proof is straightforward by the proof of Lemma \ref{corbarvkappah}.\ep\\

Form now on, we concentrate on the approximate solution $\bar v^{\kappa,h}$ which is obtained from strictly monotone scheme \ref{overlinescheme} through \reff{vbarukappah}. 
In order to provide the result, we need to use the consistency of the scheme for the regular approximate solutions. Then, the comparison principle for the scheme provides bounds over the difference between $u^{\kappa,h}$ and regular approximate solutions. 
Let
 \b*
 \Rc_{\kappa,h}[\psi](t,x)
 &:=&
 \frac{\psi(t,x)-\overline{\TT}_{\kappa,h}[\psi](t,x)}{h}
 \;+\;
 \mathcal{L}^X\psi(t,x)\\
 &&\hspace*{2.5cm}+\overline{F}_\kappa(\cdot,\psi,D\psi,D^2\psi,\psi(t,\cdot))(t,x).
 \e*
\begin{Lemma}\label{lemrateconsjump}
For a family $\{\varphi_\eps\}_{0<\eps<1}$ of smooth functions satisfying 
\be\label{regunderweps}
\left|\partial_t^{\beta_0}D^\beta\varphi^\eps\right| 
 \le C\eps^{1-2\beta_0-|\beta|_1}
 &\mbox{for any}&(\beta_0,\beta)\in\N\times\N^d\setminus\{0\},
\ee
 where $|\beta|_1:=\sum_{i=1}^d\beta_i$, and $C>0$ is some constant, we have: 
 \b*
 \left|\Rc_{\kappa,h}[\varphi_\eps]\right|_\infty
 \;\le\;
 R(h,\eps)&:=&C\;\left(h\eps^{-3}+h\theta_\kappa\eps^{-1}+h\sqrt{\theta_\kappa}+\eps^{-1}\int_{\{|z|\le\kappa\}}\!\!\!\!\!\!\!\!|z|^2\nu(dz)\right),
 \e*
for some constant $C>0$ independent of $\kappa$.
If in addition
\b*
\limsup_{h\to0}h\theta_{\kappa_h}^2<\infty~~\text{and}~~\limsup_{h\to0}\sqrt{h}\int_{\{|z|\le\kappa\}}\!\!\!\!\!\!\!\!|z|^2\nu(dz)<\infty,
\e*
we have: 
 \b*
 \left|\Rc_{\kappa_h,h}[\varphi_\eps]\right|_\infty
 \;\le\;
 R(h,\eps)&:=&C\;(h\eps^{-3}+\sqrt{h}\eps^{-1}).
 \e*
\end{Lemma}
\proof $\Rc_{\kappa,h}[\varphi_\eps]$ is bounded by 
\b*
\sup_\alpha\biggl\{\Bigl|\E\Bigl[\frac{1}{h}\bigl(\varphi_\eps(t+h,X_h^{t,x,\kappa})-\varphi_\eps(t,x)\bigr)+\frac12\Tr{a^\alpha (D^2\varphi_\eps(t+h,X_h^{t,x,\kappa})-D^2\varphi_\eps(t,x))}\\+b^\alpha (D\varphi_\eps(t+h,X_h^{t,x,\kappa})-D\varphi_\eps(t,x))+(\theta_\kappa+c^\alpha) (\varphi_\eps(t+h,X_h^{t,x,\kappa})-\varphi_\eps(t,x))\\
+\Ic^{\alpha}[\varphi_\eps](t,x)-\Ic_{\kappa,h}^{\alpha}[\varphi_\eps](t+h,x)\Bigr]\Bigr|\biggr\}
\e*
For the L\'evy integral term by Lemma \ref{lemFapproxinfinite}, we have:
\b*
|\Ic^{\alpha}[\varphi_\eps](t,x)-\Ic_{\kappa,h}^{\alpha}[\varphi_\eps](t+h,x)|
&\le&
C\biggl(|D\varphi_\eps|_\infty\bigl(\sqrt{h}+h\int_{\{|z|>\kappa\}}|z|\nu(dz)\bigr)\\
 && +h|\partial_t D^2\varphi_\eps|_\infty+|D^2\varphi_\eps|_\infty\int_{\{|z|\le\kappa\}}\!\!\!\!\!\!\!\!|z|^2\nu(dz)\biggr)\\
&\le&
C\left(h\eps^{-3}+h\sqrt{\theta_\kappa}+\eps^{-1}\int_{\{|z|\le\kappa\}}\!\!\!\!\!\!\!\!|z|^2\nu(dz)\right).
\e*
By the same argument as Lemma 3.22 in \cite{ftw} all the other terms are bounded by $h\eps^{-3}$ except
\b*
\theta_h\left(\varphi_\eps(t+h,X_h^{t,x,\kappa})-\varphi_\eps(t,x)\right)
\e*
which is bounded by $\theta_h h\eps^{-1}$.
The second assertion of the Lemma is straightforward.\ep 
\vspace*{5mm}
 
Next we need to have maximum principle for scheme \ref{overlinescheme}. Note that Lemma 3.21 in \cite{ftw} holds true for scheme \ref{overlinescheme} with $\lambda=\theta_\kappa+C$ and $\beta>\theta_\kappa+C$ where $C=\sup_\alpha|c^\alpha|$. Therefore, Proposition 3.20 in \cite{ftw} holds true for nonlocal case. More precisely, we have the following Proposition.
\begin{Proposition} \label{propmaxprinl}
Let Assumption $\FF$ holds true, and consider two arbitrary bounded functions $\varphi$ and $\psi$ satisfying:
 \b*
 h^{-1}\left(\varphi-\overline{\TT}_h[\varphi]\right) \;\le\; g_1
 &\text{and}&
 h^{-1}\left(\psi-\overline{\TT}_h[\psi]\right) \;\ge\; g_2
 \e*
for some bounded functions $g_1$ and $g_2$. Then, for every $i=0,\cdots,n$:
 \b*
 (\varphi-\psi)(t_i,x)
 &\le&
 e^{(\theta_\kappa+C_1)}|(\varphi-\psi)^+(T,\cdot)|_\infty
 +(T-h)e^{(\theta_\kappa+C_1)(T-t_i)}|(g_1-g_2)^+|_\infty~~~~~~~~~~
 \e*
where $C_1>\sup_\alpha|c^\alpha|$.
\end{Proposition}
\vspace*{5mm}

The approximation of the solution of nonlocal PDE by the Krylov method of shaking coefficients and switching system is developed in \cite{biswasjakobsenkarlsen1}. \cite{biswasjakobsenkarlsen1} provides the result of rate of convergence of general monotone schemes for the nonlocal PDE's satisfying Assumption IHJB. However, they referred the regularity of the approximate solutions to the result of \cite{biswasjakobsenkarlsen2} where the approximate solution obtained from switching system could only be locally $\frac12$-H\"older continuous on $t$. But, in the case of scheme \reff{schemejump}, we need the solution of \reff{equationnl}-\reff{terminalnl} be uniformly $\frac12$-H\"older continuous on $t$. It is because we need the regular approximate solutions obtained from Krylov method and switching solution to satisfy \reff{regunderweps}. Therefore, in the present work we need to rebuild Lemma 5.3 in \cite{biswasjakobsenkarlsen2} under the Assumption IHJB to obtain global $\frac12$-H\"older continuous on $t$ for the solution of the switching system.

Therefore, we continue this subsection by introducing the switching system of nonlocal PDE's with the regularity result needed for the solution of this system.\\
Let $k$ be a non-negative constant. Suppose the following system of PDE's:
\be\label{switchsys}
 \max\left\{-\Lc^X v_i(t,x)- F_i\left(t,x,v_i(t,x),D v_i(t,x),D^2v_i(t,x),v_i(t,\cdot)\right),v_i-\Mc^iv\right\}=0\nonumber\\
 v_i(T,\cdot)=g_i(\cdot),~~~~~~~~~~~~~
\ee
where $i=1,\cdots,M$ and
\b*
 F_i(t,x,r,p,\gamma,\psi)
 &:=&
 \mathop{\inf}\limits_{\alpha\in\Ac_i}\left\{\Lc^\alpha(t,x,r,p,\gamma,\gamma)+\Ic^\alpha(t,x,r,p,\gamma,\psi)\right\}\\
 \Lc^\alpha(t,x,r,p,\gamma,\gamma)
 &:=&
 \frac{1}{2}\Tr{a^\alpha(t,x)\gamma}+b^\alpha(t,x)\cdot p+c^\alpha(t,x)r+k^\alpha(t,x)\\
 \Ic^\alpha(t,x,r,p,\gamma,\psi)
 &:=&
 \int_{\R^d_*}\left(\psi\left(t,x+\eta^\alpha(t,x,z)\right)-r-\mathds{1}_{\{|z|\le1\}}\eta^\alpha(t,x,z)\cdot p\right)d\nu(z)\\
 \Mc^ir
 &:=&
 \mathop{\min}\limits_{j\ne i}r_j+k.
\e*
We would like to emphasize that $g_i$s need to satisfy $g_i-\Mc^i\bar g\le0$ where $\bar g=(g_1,\cdots,g_M)$. If each $g_i=g$ then we obviously have $g_i-\Mc^i\bar g\le0$.

Existence and comparison principle result for the above switching system is provided in Proposition 6.1 \cite{biswasjakobsenkarlsen1}. Also, it is known from Theorem 6.3 in \cite{biswasjakobsenkarlsen1}, that if $(v^1,\cdots,v^M)$ and $v$ be  respectively the solutions of \reff{switchsys} and \reff{equationnl}-\reff{terminalnl} with $\Ac=\cup_{i=1}^M\Ac_i$ and $\Ac_i$s are disjoint sets, then 
\be 
\label{switchapporx}
0\le v^i-v\le Ck^\frac{1}{3}~~\text{for}~~i=1,\cdots,M.
\ee
The regularity result for \reff{switchsys} is provided in \cite{biswasjakobsenkarlsen2}. There, it is proved that $(v^i)$ is Lipschitz with respect to $x$ and locally 1/2-H\"older continuous with respect to $t$. For the proof of Theorem \ref{thmrateconvjump}, $(v^i)$ should be uniformly 1/2-H\"older continuous with respect to $t$. The following Lemma provide the uniform 1/2-H\"older continuity for $(v^i)$.
\begin{Lemma}\label{lemliphold}
Assume {\bf HJB} holds for each $i$ and let $(v^i)$ be the viscosity solution of \reff{switchsys}. Then there exist a constant $C$ such that for any $i=1,\cdots,M$:
\b*
\left|v^i\right|_1&\le& C.
\e* 
\end{Lemma}
\proof Lipschitz continuity with respect to $x$ is done in Lemma 5.2 in \cite{biswasjakobsenkarlsen2}. To obtain uniform $1/2-$H\"older continuity with respect to $t$, the proof of Lemma 5.3 in \cite{biswasjakobsenkarlsen2} should be modified  by using assumption {\bf HJB}.\\
Fix $y\in\R^d$, $t$ and $t'$ where $t\le t'$. For each $i=1,\cdots,M$, define:
\b*
\psi_i(t,x)&:=&\lambda \frac{L}{2}\left[e^{A(t'-t)}|x-y|^2+B(t'-t)\right] + K(t'-t) + \lambda^{-1}\frac{L}{2} + v^i(t',y)
\e*
Where $L=\frac{1}{2}|v|_1$ and $\lambda$, $a$ and $\gamma$ will be defined later. Then:
\b*
 \partial_t \psi_i(t,x)
 &=&
 -\lambda \frac{L}{2}\left(Ae^{A(t'-t)}|x-y|^2+B\right)-K
 \\
 D\psi_i(t,x)
 &=&
 \lambda Le^{A(t'-t)}(x-y)
 \\
 D^2\psi(t,x)
 &=&
 \lambda Le^{A(t'-t)}I_{d\times d}.
\e* 
So, 
\b*
\lefteqn{-\partial_t \psi_i +\Lc^\alpha(t,x,\psi_i,D\psi_i,D^2\psi_i)+\Ic^\alpha(t,x,\psi_i,D\psi_i)}\\
&=& -\lambda L\left(Ae^{A(t'-t)}|x-y|^2+B\right)-K\\
&&+\lambda Le^{A(t'-t)}\Tr{a^\alpha(t,x)}+2\lambda Le^{A(t'-t)}b^\alpha(t,x)\cdot(x-y)+c^\alpha(t,x)\psi_i+k^\alpha(t,x)\\
&&+ \lambda \frac{L}{2}e^{A(t'-t)}\int_{\R^d_*}\left(|x+\eta^\alpha(t,x,z)-y|^2-|x-y|^2-21\!\!1_{\{|z|\le1\}}\eta^\alpha(t,x,z)\cdot (x-y)\right)d\nu(z).
\e*
By {\bf HJB}, we choose $K$ and $\lambda$ so that,
\b*
|a^\alpha|_\infty\le K, |b^\alpha|_\infty\le K, |c^\alpha|_\infty\le K, |k^\alpha|_\infty\le K, K^{-1}\le\lambda\le K\\
|v|_\infty\le K, |\eta^\alpha(t,x,z)|\le K(1\wedge|z|).
\e*
Without loss of generality and with the similar argument as in Remark 3.7, we  can suppose that for any $\alpha$, $c^\alpha\le0$. So by choosing positive large $D$, there exists positive constants $C_1$, $C_2$ and $C_3$ such that:
\b*
\lefteqn{-\partial_t \psi_i +\Lc^\alpha(t,x,\psi_i,D\psi_i,D^2\psi_i)+\Ic^\alpha(t,x,\psi_i,D\psi_i)}\\
&\le& -\lambda Le^{A(t'-t)}K\left(\left(\frac{A}{K}-\frac{1}{2}\right)|x-y|^2-\frac{3}{2}|x-y|+C_1B-C_2\right)+C_3.
\e*
Therefore, choice of large $B$ and $D$ makes the right hand side negative.
\b*
-\partial_t \psi_i +\Lc^\alpha(t,x,\psi_i,D\psi_i,D^2\psi_i)+\Ic^\alpha(t,x,\psi_i,D\psi_i)
\le 0.
\e*
On the other hand,
\b*
\psi(t',x)&=&\frac{L}{2}\left(\lambda|x-y|^2+\lambda^{-1}\right)+v^i(t',y).
\e*
Minimizing with respect to $\lambda$,
\b*
\psi(t',x)\ge L|x-y|+v^i(t',y)\ge v^i(t',x).
\e*
We can conclude that $\psi_i$ is a super solution of \reff{switchsys}. So, by comparison Theorem in \cite{biswasjakobsenkarlsen1},
\b*
\psi_i(t,y)\ge v^i(t,y).
\e*
So,
\b*
\frac{L}{2}\left(\lambda B(t'-t)+\lambda^{-1}\right)+v^i(t',y)\ge v^i(t,y).
\e*
Therefore,
\b*
v^i(t,y)-v^i(t',y)\le C\sqrt{t'-t}.
\e*
The other inequality can be done similarly by choosing:
\b*
\psi_i(t,x)&:=&-\lambda \frac{L}{2}\left[e^{A(t'-t)}|x-y|^2-B(t'-t)\right] - K(t'-t) - \lambda^{-1}\frac{L}{2} + v^i(t',y).
\e*
\ep\\
\begin{Remark}
{\rm
Notice that all the result of switching system is correct for \reff{equationnl}-\reff{terminalnl} satisfying {\bf F} by simply setting $M=1$ and $k=0$.
}
\end{Remark}

Therefore, by \cite{biswasjakobsenkarlsen1} there are regular functions $\underline{w}^\kappa_\eps$ and $\overline{w}^\kappa_\eps$ which are respectively the regular sub- and super-solution of 
\b*
 &&-\Lc^X u^\kappa(t,x)- \overline{F}_\kappa\left(t,x,u^\kappa(t,x),D u^\kappa(t,x),D^2u^\kappa(t,x),u^\kappa(t,\cdot)\right)= 0,~~\mbox{on}~[0,T)\times\R^d,~~~~~~~~~~~~~~~~~~~~~~\\
&&u^\kappa(T,\cdot)=g, ~~~~~~~~~~~~~~~~~~~~~~~~~~~~~~~~~~~~~~~~~~~~~~~~~~~~~~~~~~~~~~~~~\mbox{on}~\in\R^d. \nonumber
\e*
where
\b*
 \overline{F}_\kappa(t,x,r,p,\gamma,\psi)
 &:=&
\inf_{\alpha\in\Ac}\left\{\Lc^{\alpha,\beta}(t,x,r,p,\gamma)+\Ic_\kappa^{\alpha,\beta}(t,x,r,p,\gamma,\psi)\right\}
 \e*
 (one can replace $\sup\inf$ by $\inf\sup$) where
   \b*
 \Lc^{\alpha,\beta}(t,x,r,p,\gamma)
 &:=&
 \frac{1}{2}
 \Tr{\sigma^{\alpha,\beta}\sigma^{{\alpha,\beta}{\rm T}}(t,x)\gamma}
 +b^{\alpha,\beta}(t,x)p+(c^{\alpha,\beta}(t,x)+\theta_\kappa)r,
 \e*
and
 \b*
 \Ic_\kappa^{\alpha,\beta}(t,x,r,p,\gamma,\psi)&\!\!\!\!\!:=\!\!\!\!\!&\int_{\{|z|>\kappa\}}\!\!\!\!\!\!\!\!\!\bigl(\psi\bigl(x+\eta^{\alpha,\beta}(t,x,z)\bigr)-r-\mathds{1}_{\{|z|\le1\}}\eta^{\alpha,\beta}(t,x,z)\cdot p\bigr)\nu(dz).
 \e*
Therefore, by Proposition 6.2 and Theorem 6.3  of \cite{biswasjakobsenkarlsen1}, Lemma \ref{lemrateconsjump} and Proposition \ref{propmaxprinl},
\b*
(u^\kappa-u^{\kappa,h})(t,x)&\le&(u^\kappa-\underline{w}^\kappa_\eps+\underline{w}^\kappa_\eps-u^{\kappa,h})(t,x)\\
&\le&
Ce^{(\theta_\kappa+C_1)(T-t)}\left(\eps+h\eps^{-3}+h\theta_\kappa\eps^{-1}+h\sqrt{\theta_\kappa}+\eps^{-1}\int_{\{|z|\le\kappa\}}\!\!\!\!\!\!\!\!|z|^2\nu(dz)\right)
\e*
and
\b*
(u^{\kappa,h}-u^\kappa)(t,x)&\le&(u^{\kappa,h}-\overline{w}^\kappa_\eps+\overline{w}^\kappa_\eps-u^\kappa)(t,x)\\
&\le&
Ce^{(\theta_\kappa+C_1)(T-t)}\left(\eps^\frac{1}{3}+h\eps^{-3}+h\theta_\kappa\eps^{-1}+h\sqrt{\theta_\kappa}+\eps^{-1}\int_{\{|z|\le\kappa\}}\!\!\!\!\!\!\!\!|z|^2\nu(dz)\right).
\e*
Notice that $v^\kappa(t,x)=e^{-\theta_\kappa(T-t)}u^\kappa(t,x)$. So,
\b*
v^\kappa-\bar v^{\kappa,h}&\le&C\left(\eps+h\eps^{-3}+h\theta_\kappa\eps^{-1}+h\sqrt{\theta_\kappa}+\eps^{-1}\int_{\{|z|\le\kappa\}}\!\!\!\!\!\!\!\!|z|^2\nu(dz)\right)
\e*
and
\b*
\bar v^{\kappa,h}-v^\kappa&\le&C\left(\eps^\frac{1}{3}+h\eps^{-3}+h\theta_\kappa\eps^{-1}+h\sqrt{\theta_\kappa}+\eps^{-1}\int_{\{|z|\le\kappa\}}\!\!\!\!\!\!\!\!|z|^2\nu(dz)\right).
\e*
On the other hand, because of \reff{kappahrate} and by Lemma \reff{corbarvkappah}, the second part of Theorem \ref{thmrateconvjump} is provided after choice of optimal $\eps$.

\section{Conclusion}
The algorithm is the first probabilistic numerical method for fully nonlinear nonlocal parabolic problems. As in local case (\cite{ftw}), it converges to the viscosity solution of the problem. A rate of convergence is known for the convex (concave) nonlinearities. Also with the same argument as in Section 4 in \cite{ftw}, Monte Carlo approximations of expectations inside the scheme do not affect the asymptotic results if enough number of samples would be used. The error analysis for MCQ shows that the appropriate approximation of jump-diffusion process with compound Poisson process could be applied in discretization procedure. The theoretical result is followed by some numerical examples which confirms the convergence of the scheme.

On the other hand there are some features where the scheme is not implementable in nonlocal case, e.g. when the nonlinearity is of HJB type. This could be the challenge of future works. The other open issue is to relax some assumptions. For example, relaxing the assumption of uniform ellipticity may be a future work.

\bibliographystyle{plain}
\bibliography{bib}

\end{document}